\documentclass{conm-p-l}
\usepackage{amsmath, amstext, amsbsy, amssymb, amscd, graphics}

\newtheorem{theorem}{Theorem}[section]
\newtheorem{lemma}[theorem]{Lemma}
\newtheorem{proposition}[theorem]{Proposition}

\theoremstyle{definition}

\theoremstyle{remark}
\newtheorem{remark}[theorem]{Remark}

\numberwithin{equation}{section}

\newcommand{\C}{ \mathbb C }
\newcommand{\ch}{{\rm ch}}
\newcommand{\duval}{{\C^2//\G}}
\newcommand{\End}{{\rm End}}
\newcommand{\FFock}{{\mathcal F}}
\newcommand{\g}{{\gamma}}
\newcommand{\G}{{\Gamma}}
\newcommand{\Gn}{{\Gamma}_n}

\newcommand{\Hn}{{\mathbb H}_n}
\newcommand{\Hx}{{\mathbb H}}
\newcommand{\la}{\lambda}

\newcommand{\Pnr}{\mathcal P_n (r)}

\newcommand{\RG}{ \mathbb R_\G}
\newcommand{\Supp}{\text{\rm Supp}}

\newcommand{\td}{{\rm td}}

\newcommand{\vac}{|0\rangle}
\newcommand{\Xn}{ X^{[n]} }
\newcommand{\Z}{ \mathbb Z}

{\vskip-\lastskip\medskip
  \noindent
  {\em #1.}\enspace
  }%
{\qed\par\medskip
  }

\begin{document}
\title[Hilbert schemes and soliton equations]
  {Hilbert schemes of points on the minimal resolution and soliton
equations}

\author[Zhenbo Qin]{Zhenbo Qin}
\address{Department of Mathematics, University of Missouri, Columbia,
MO 65211, USA} \email{zq@math.missouri.edu}

\author[Weiqiang Wang]{Weiqiang Wang}
\address{Department of Mathematics, University of Virginia,
Charlottesville, VA 22904} \email{ww9c@virginia.edu}
\thanks{Both authors are partially supported by NSF grants}


\begin{abstract}
The equivariant and ordinary cohomology rings of Hilbert schemes
of points on the minimal resolution $\duval$ for cyclic $\G$ are
studied using vertex operator technique, and connections between
these rings and the class algebras of wreath products are
explicitly established. We further show that certain generating
functions of equivariant intersection numbers on the Hilbert
schemes and related moduli spaces of sheaves on $\duval$ are
$\tau$-functions of $2$-Toda hierarchies.
\end{abstract}

\maketitle
\date{}

\section{\bf Introduction}

In the past few years, the theory of vertex operators (cf.
\cite{FLM}) has found remarkable applications in the study of the
Hilbert schemes $\Xn$ of points on a surface $X$ (see \cite{Na2,
Le, LQW1, LQW2, Ru, Vas, Wa1} and the references therein). To a
large extent, the construction of Heisenberg algebra by Nakajima
\cite{Na1} (also cf. \cite{Gro}) was a key starting point. It is
well known that the Hilbert-Chow morphism from the Hilbert scheme
$\Xn$ to the symmetric product $X^n/S_n$ is a resolution of
singularities. An isomorphism has been established between the
cohomology rings of the Hilbert schemes and the Chen-Ruan orbifold
cohomology rings of the symmetric products for the affine plane
$\C^2$ in \cite{LS, Vas} via different approaches, and for a large
class of quasi-projective surfaces in \cite{LQW1}. These and other
results have supported Ruan's general conjectures \cite{Ru} on the
relations between the orbifold cohomology ring of an orbifold and
the cohomology ring of its crepant resolution. On the other hand,
a certain generating function of equivariant intersection numbers
on the Hilbert schemes $(\C^2)^{[n]}$ are shown \cite{LQW2} to be
$\tau$-functions of the $2$-Toda hierarchies.

When $X =\duval$ is the minimal resolution of $\C^2/\G$, where
$\G$ is a finite subgroup of $SL_2(\C)$, we have the following
resolution of singularities \cite{Wa1}
$$\pi_n: (\duval)^{[n]} \longrightarrow \C^{2n}/ \Gn$$
obtained by the composition $(\duval)^{[n]} \rightarrow
(\duval)^{n}/ S_n \rightarrow (\C^{2}/\G)^n/ S_n \cong \C^{2n}/
\Gn$.
Here $\Gn := \G^n \rtimes S_n$ is the wreath product. When $\G$ is
cyclic, we fix a distinguished $T=\C^*$ action on $\duval$ which
induces an action on $(\duval)^{[n]}$ with isolated fixed points.
Throughout the paper, we will assume that $\G$ is a cyclic finite
subgroup of $SL_2(\C)$. For example, if $\G$ is of order $2$, then
$\duval$ is isomorphic to the total space of the cotangent bundle
over the projective line $\mathbb P^1$.

A main goal of this paper is to study the (equivariant) cohomology
ring of Hilbert schemes of points on $\duval$ using vertex
operator technique and to develop its relations to soliton
equations. This study specializes when $\G$ is trivial to the
study of \cite{Vas, LQW2} in the affine plane case. In addition,
we shall establish Ruan's conjecture on the cohomology ring
isomorphism for the resolution $\pi_n: (\duval)^{[n]} \to \C^{2n}/
\Gn$ and provide an explicit map for such an isomorphism.

Let us now discuss the paper in more details. First of all,
generalizing the affine plane case studied by Vasserot \cite{Vas},
we introduce a ring structure on
\begin{eqnarray*}
\Hn =H^{2n}_T ((\duval)^{[n]})
\end{eqnarray*}
which encodes the equivariant cohomology ring structure of $H^*_T
((\duval)^{[n]})$. We further construct a Heisenberg algebra
acting irreducibly on
\begin{eqnarray*}
\mathbb H = \bigoplus_{n = 0}^\infty \Hn.
\end{eqnarray*}
This is an equivariant analog of the Heisenberg algebra in
\cite{Na1, Na2} (also cf. \cite{Vas}). Our study uses in
an essential way a very concrete and useful description of $\duval$
and its torus action provided by Ito-Nakamura \cite{IN}.

Next, we introduce an explicit map from the ring $\Hn$ to the
class algebra $R(\Gn)$ and show that it is a ring isomorphism.
This result specializes to a result in \cite{Vas} in the affine
plane case. Denote by $\mathcal G_\G^*(n)$ the graded ring
associated to a natural filtration on $R(\Gn)$ (cf. \cite{Wa3}).
We obtain an explicit map from $H^*( (\duval)^{[n]})$ to $\mathcal
G_\G^*(n)$ which is further shown to be also a graded ring
isomorphism. When $\G$ is trivial, this has been established in
\cite{LS, Vas} and in \cite{LQW1} using different methods. With
the help of the Heisenberg algebra in the setup of wreath products
\cite{FJW, Wa2} and the Heisenberg algebra on Hilbert schemes
constructed above, our maps further identify several distinguished
linear bases on both sides. Our above isomorphism of graded rings
establishes Ruan's conjecture for the crepant resolution $\pi_n:
(\duval)^{[n]} \rightarrow \C^{2n}/ \Gn$, since one can identify
$\mathcal G_\G^*(n)$ as the Chen-Ruan orbifold cohomology ring of
$\C^2/\Gn$. The resolution $\pi_n$ seems to be the first
nontrivial example beyond the Hilbert-Chow morphism where Ruan's
conjecture is established in a constructive way (see \cite{EG} for
an earlier nonconstructive approach toward Ruan's conjecture). A
direct and constructive ring isomorphism from $H^*(
(\duval)^{[n]})$ to $\mathcal G_\G^*(n)$ for non-cyclic $\G$ has
yet to be constructed, though $\mathcal G_\G^*(n)$ can still be
defined for a general $\G$, cf. \cite{EG, Wa3}.

In addition, we introduce a family of moduli spaces of sheaves on
$\duval$, which are isomorphic to the Hilbert schemes
$(\duval)^{[n]}$ and parameterized by certain integral lattice in
$H^2_T(\duval)$. Using an operator approach, we study the
$T$-equivariant Chern characters of some distinguished
$T$-equivariant tautological bundles over these moduli spaces. We
identify the Chern character operators, which are defined in terms
of the $T$-equivariant Chern characters, with some familiar
operators acting on the fermionic Fock space associated to the
integral lattice in $H^2_T(\duval)$. We then formulate generating
functions of the equivariant intersection numbers of these Chern
characters, and recast them in an operator formalism. It follows
from standard arguments that they are $\tau$-functions of the
$2$-Toda hierarchies of Ueno-Takasaki \cite{UT}. These results
generalize the affine plane case studied in \cite{LQW2}.

This paper is organized as follows. In Section~\ref{sec:hilb}, we
study the torus-equivariant geometry of the Hilbert schemes and
formulate the Heisenberg algebra. In Section~\ref{sec:isom}, we
established the two ring isomorphisms in the equivariant and
ordinary cohomology setups. In Section~\ref{sec:moduli}, we
introduce a certain moduli space of sheaves on the minimal
resolution $\duval$ and study the related fermionic Fock space. In
Section~\ref{sec:inter}, we show that certain generating functions
of the equivariant intersection numbers on these moduli spaces are
$\tau$-functions.

\section{\bf Hilbert schemes of points on the minimal resolution}
\label{sec:hilb}
\subsection{The minimal resolution}

Let $\G =\langle a \rangle$ be the cyclic group of order $r$.
Below we often regard $\G$ as the subgroup of $SL_2(\C)$ which
consists of the diagonal matrices $a^i = \text{diag}
(\varepsilon^i,\varepsilon^{-i}), 0 \le i <r$, where $\varepsilon
=e^{2\pi \sqrt{-1}/r}$ is a primitive $r$-th root of unity. We
denote by $\Gamma^*= \{\g_0,\g_1, \ldots, \g_{r-1}\}$ the set of
complex irreducible characters of $\G$, with the character table
of $\G$ given by
$$ \g_k (a^i) =\varepsilon^{ik}, \qquad 0 \le i,k \le r-1. $$
In other words, the character table is given by the $r \times r$
matrix
$$ M = [ \varepsilon^{ik} ]_{0 \le i,k \leq r-1}. $$
Let $\pi: \duval \rightarrow \C^2 /\G$ be the minimal
desingularization of $\C^2 /\G$, and let $X =\duval$. The
exceptional fiber $\pi^{-1} (0)$ consists of $(r-1)$ projective
lines $\Sigma_1, \ldots, \Sigma_{r-1}$, with the configuration of
a Dynkin diagram
%
 \begin{equation*}
 \begin{picture}(50,35) 
 \put(-90,18){$A_{r-1}:$}
 \put(-22,20){$\circ$}
 \put(-16,23){\line(1,0){18}}
 \put(3,20){$\circ$}
 \put(9,23){\line(1,0){18}}
 \put(28,20){$\circ$}
  \put(34,22){ \dots }
 \put(64,23){\line(1,0){18}}
 \put(83,20){$\circ$}
 \put(-25,5){$\Sigma_1$}
 \put(2,5){$\Sigma_2$}
 \put(80,5){$\Sigma_{r-1}$}
 \end{picture}
 \end{equation*}
in the following sense: $\Sigma_i$ and $\Sigma_j$ intersect if and
only if $|i-j|=1$; and when this is so they intersect
transversally at one point. Following Ito-Nakamura \cite{IN}, we
can identify $X$ with the subvariety of $(\C^2)^{[r]}$ which
consists of the points corresponding to the following
$\G$-invariant ideals in the coordinate ring $\C[z_1, z_2]$ of
$\C^2$:
\begin{eqnarray}
I(O) &:=&\prod_{(a_1, a_2) \in O} \mathfrak m_{(a_1, a_2)}
    = (z_1^r - a_1^r, z_1z_2 - a_1a_2, z_2^r - a_2^r),
    \label{ideals1}   \\
I_i(p_i: q_i) &:=&(p_iz_1^i - q_iz_2^{r-i}, z_1z_2, z_1^{i+1},
    z_2^{r+1-i})   \label{ideals2}
\end{eqnarray}
where $\mathfrak m_{(a_1, a_2)} = (z_1 - a_1, z_2 - a_2)$, $O$
stands for a $\G$-orbit in $\C^2$ disjoint from the origin, $1 \le
i \le r-1$, and $(p_i : q_i) \in \mathbb P^1$. For each $1 \le i
\le r-1$, the rational curve $\Sigma_i$ consists of all the points
in $(\C^2)^{[r]}$ corresponding to the ideals (\ref{ideals2}).

Moreover, $X$ admits an affine open cover $\{X_i\}_{0 \le i \le
r-1}:$
\begin{eqnarray*}
 X =\cup_{0 \le i \le r-1} X_i, \quad \text{where } X_i \cong
\text{\rm Spec} \,\, \C[z_{i, 1}, z_{i, 2}].
\end{eqnarray*}
The inclusion map $X_i \hookrightarrow X$ is given by the morphism
defined via the universal property of $(\C^2)^{[r]}$ from the
following $2$-dimensional flat family of $\G$-invariant ideals in
$\C[z_{1}, z_{2}]$:
\begin{eqnarray}  \label{ideal_X_i}
\mathfrak I_i(a_{i, 1}, a_{i, 2}) := (z_1^{i+1} - a_{i,
1}z_2^{r-i-1}, z_1z_2 - a_{i, 1}a_{i, 2}, z_2^{r-i} - a_{i,
2}z_1^{i})
\end{eqnarray}
where $(a_{i, 1}, a_{i, 2})$ stands for points in $X_i \cong
\C^2$. Note that when $a_{i, 1}a_{i, 2} \ne 0$, we have
$\mathfrak I_i(a_{i, 1}, a_{i, 2}) = I(O)$ where $O$ is the
$\G$-orbit of the point $(a_1, a_2) \in \C^2$ with
\begin{eqnarray*}
a_1 = a_{i, 1}^{(r-i)/r}a_{i, 2}^{(r-i-1)/r}, \quad \quad
a_2 = a_{i, 1}a_{i, 2}/a_1
 =a_{i, 1}^{i/r}a_{i, 2}^{(i+1)/r}
\end{eqnarray*}
(i.e., $a_{i, 1} = a_1^{i+1}/a_2^{r-i-1}$ and $a_{i, 2} =
a_2^{r-i}/a_1^{i}$). Similarly, we obtain
\begin{eqnarray}
\mathfrak I_i(a_{i, 1}, 0)
 &=& I_{i+1}(1: a_{i, 1}), \quad\quad 0 \le i \le r-2
 \label{axis1}  \\
\mathfrak I_{r-1}(a_{r-1, 1}, 0)
 &=& (z_1^r - a_{r-1, 1}, z_2),
 \label{axis2}  \\
\mathfrak I_i(0, a_{i, 2})
 &=& I_{i}(a_{i, 2}: 1),\qquad\quad  1 \le i \le r-1
 \label{axis3}  \\
\mathfrak I_0 (0, a_{0, 2})
 &=& (z_1, z_2^r - a_{0, 2}).   \label{axis4}
\end{eqnarray}

Let $\xi_0, \ldots, \xi_{r-1} \in X$ be the points corresponding
respectively to the ideals:
\begin{align}  \label{fixed_pts}
I_1(1: 0), \,\, I_1(0: 1) = I_2(1: 0), \ldots, \,\, I_{r-2}(0: 1)
= I_{r-1}(1: 0), \,\, I_{r-1}(0: 1).
\end{align}
Then $\xi_0 \in \Sigma_1$, $\{ \xi_1 \} = \Sigma_1 \cap \Sigma_2,
\ldots, \{ \xi_{r-2} \} = \Sigma_{r-2} \cap \Sigma_{r-1}$, and
$\xi_{r-1} \in \Sigma_{r-1}$. Moreover, we see from
(\ref{axis1})--(\ref{axis4}) that $\xi_i$ is the origin of the
open affine chart $X_i$, $\Sigma_{i+1} - \{ \xi_{i+1} \}$ is the
$z_{i,1}$-axis of $X_i$ when $0 \le i \le r-2$, and $\Sigma_{i} -
\{ \xi_{i-1} \}$ is the $z_{i,2}$-axis of $X_i$ when $1 \le i \le
r-1$. Let $\Sigma_0$ be the $z_{0,2}$-axis of $X_0$, and
$\Sigma_r$ be the $z_{r-1, 1}$-axis of $X_{r-1}$.

\subsection{Fixed point classes and bilinear forms}

Let $T=\C^*$ act on $\C^2$ by
\begin{eqnarray}     \label{T_action_C^2}
s(z_1, z_2) = (s z_1, s^{-1}z_2), \quad s \in T.
\end{eqnarray}
This action induces a $T$-action on $(\C^2)^{[r]}$ which preserves
the minimal resolution $X = \duval$ as a subvariety of
$(\C^2)^{[r]}$. In view of (\ref{ideals1}) and (\ref{ideals2}),
the locus $X^T$ of $T$-fixed points is
\begin{eqnarray*}
X^T = \{ \xi_0, \ldots, \xi_{r-1} \}.
\end{eqnarray*}

Since the family of ideals in (\ref{ideal_X_i}) is $T$-invariant,
$X_i$ is $T$-invariant. In addition, $T$ acts on points $(a_{i,
1}, a_{i, 2}) \in X_i$ by $s(a_{i, 1}, a_{i, 2}) = (s^{-r}a_{i,
1}, s^r a_{i, 2})$, i.e., $T$ acts on the coordinate functions
$(z_{i, 1}, z_{i, 2})$ of $X_i$ by
\begin{eqnarray}  \label {action_on_X_i}
s(z_{i, 1}, z_{i, 2}) = (s^{r}z_{i, 1}, s^{-r} z_{i, 2}).
\end{eqnarray}

The $T$-action on $X$ induces a $T$-action on the Hilbert scheme
$\Xn$, which again has isolated fixed points. As explained below,
these isolated fixed points $\xi_\la$ are parametrized by the
multi-partitions $\la$ in the finite set
\begin{eqnarray}  \label{Pnr}
\Pnr = \{\la=(\la^0, \cdots, \la^{r-1}) \mid  |\la^0|+ \cdots +
|\la^{r-1}| =n \}.
\end{eqnarray}
More explicitly, a $T$-fixed subscheme $Z$ of $X$ of length-$n$ is
a disjoint union of $T$-fixed subschemes $Z_i$, $i=0, \dots, r-1$,
such that $Z_i$ is supported at $\xi_i$, and $\sum_{i} \ell(Z_i) =
n$. Put $n_i = \ell(\xi_i)$. Since $\xi_i \in X_i \cong \C^2$ is
the origin and $T$ acts on $X_i$ by (\ref{action_on_X_i}), we see
from \cite{ES} that the $T$-fixed subschemes $Z_i \in
(X_i)^{[n_i]} \subset X^{[n_i]}$ are in one-to-one correspondence
with the partitions $\la^i$ of $n_i$. We denote such $Z_i$ by
$\xi_i^{\la^i}$ and so $Z$ is given by $\xi_\la = \sum_{i}
\xi_i^{\la^i}$ when $\la = (\la^0, \cdots, \la^{r-1})$. It is
known \cite{ES} that the tangent space
$T_{\xi_i^{\la^i}}X^{[n_i]}$ of $X^{[n_i]}$ at $\xi_i^{\la^i}$ as
a $T$-module decomposes as
\begin{eqnarray*}
T_{\xi_i^{\la^i}}X^{[n_i]} = \bigoplus_{\square \in \la^i}
\left(\theta^{rh(\square)} \bigoplus \theta^{-rh(\square)} \right)
\end{eqnarray*}
where $\theta$ is the $1$-dimensional standard module of $T$, and
$\square$ runs over the cells in the Young diagram corresponding
to the partition $\la^i$, $h(\square)$ is the hook number of a
cell $\square$. Hence as a $T$-module the tangent space of
$X^{[n]}$ at ${\xi_{\la}}$ decomposes as
\begin{eqnarray}  \label{tang_space}
T_{\xi_{\la}}X^{[n]} = \bigoplus_{i=0}^{r-1} \bigoplus_{\square\in
\la^i} \left(\theta^{rh(\square)} \bigoplus \theta^{-rh(\square)}
\right).
\end{eqnarray}

Let $H^*_{T}(M)$ be the equivariant cohomology of a smooth variety
$M$ with $\mathbb C$-coefficient. Then $H^*_{T}(M)$ is a $\mathbb
C[t]$-module if we identify $H^*_{T}(\text{\rm pt})$ and $\mathbb
C[t]$, where $t$ is an element of degree-$2$.
Putting $h(\la)=\prod_{i=0}^{r-1} \prod_{\square \in \la^i}
h(\square)$, we see from (\ref{tang_space}) that the equivariant
Euler class of the tangent bundle at $\xi_{\la} \in X^{[n]}$ is
\begin{align} \label{normal}
e_T(T_{\xi_{\la}} \Xn) = (-1)^n r^{2n} t^{2n} \prod_{i=0}^{r-1}
\,\, \prod_{\square \in \la^i}h(\square)^2 = (-1)^n r^{2n} t^{2n}
h(\la)^2.
\end{align}

For $\la \in \Pnr$, let $i_{\la}: \xi_{\la} \to X^{[n]}$ be the
inclusion map. Let $[\xi_{\la}] \in H^{4n}_T(\Xn)$ be the
equivariant cohomology class corresponding to the fixed point
$\xi_{\la}$. Then, $[\xi_{\la}] = i_{\la}^!(1_{\xi_{\la}})$ where
$1_{\xi_{\la}}$ is the unity in the ring $H^*_T(\xi_{\la})$ and
$i_{\la}^!$ is the Gysin map. For $\la$ and $\mu$ in $\Pnr$, we
see from the projection formula and (\ref{normal}) that
\begin{align} \label{xi_cup}
[\xi_{\la}] \cup [\xi_{\mu}] = \delta_{\la, \mu} e_T(T_{\xi_{\la}}
\Xn) [\xi_{\la}] = \delta_{\la, \mu} \cdot (-1)^n r^{2n} t^{2n}
h(\la)^2 [\xi_{\la}].
\end{align}

Denote $\iota_n = \bigoplus\limits_{\la \in \Pnr} i_{\la}:
(X^{[n]})^T \to X^{[n]}$, and let $\iota_n^!: H^*_T((X^{[n]})^T)'
\to H^*_T(X^{[n]} )'$ be the induced Gysin map where $H^*_T( \cdot
)' = H^*_T( \cdot ) \otimes_{\mathbb C[t]} \mathbb C(t)$. By the
localization theorem, $\iota_n^!$ is an isomorphism.
The inverse $(\iota_n^!)^{-1}$ is given by
\begin{eqnarray} \label{inverse}
\alpha \to \left ( \frac{i_{\la}^*\alpha} {e_T(T_{\xi_{\la}} \Xn)}
\right )_{\la \in \Pnr}
= \left ( \frac{i_{\la}^*\alpha} {(-1)^n r^{2n} t^{2n} h(\la)^2}
\right )_{\la \in \Pnr}.
\end{eqnarray}
We define a bilinear form $\langle -, - \rangle$ on
$H^*_T(X^{[n]})' \otimes_{\mathbb C(t)}H^*_T(X^{[n]})'\to \mathbb
C(t)$ by
\begin{eqnarray}  \label{form}
\langle \alpha, \beta\rangle =(-1)^n
p_n^!(\iota_n^!)^{-1}(\alpha\cup\beta)
\end{eqnarray}
where $p_n$ is the projection of the set $(X^{[n]})^T$ of
$T$-fixed points to a point.

Set $\Hn = H^{2n}_T (\Xn)$. From the spectral sequence computation
and the fact that $H^{2k}(X^{[n]})=0$ for $k > n$, we see that
$H^{4n}_T(X^{[n]}) = t^n \cdot H^{2n}_T(X^{n]}) = t^n \Hn$.
So there is an induced (commutative
associative) ring structure $\star$ on $\Hn$ defined by:
\begin{eqnarray*}
t^n (x \star y) =x \cup y, \qquad x,y \in \Hn.
\end{eqnarray*}

Associated to $[\xi_{\la}] \in H^{4n}_T(\Xn)$, we can define
$[\la] \in \Hn$ for $\la \in \Pnr$ by
\begin{eqnarray}  \label{rescaling}
[\la]=(-1)^n r^{-n} h(\la)^{-1} t^{-n} [\xi_{\la}] \in \Hn.
\end{eqnarray}
To emphasis the case $n=1$, we introduce for $0 \le i \le r-1$ the
notation:
\begin{eqnarray*}
\diamondsuit_{i} = -r^{-1}t^{-1}[\xi_i].
\end{eqnarray*}
It follows from (\ref{xi_cup}) and (\ref{form}) that for $\la, \mu
\in \Pnr$, we have
\begin{eqnarray}  \label{pairing_Hn}
\langle [\la], [\mu] \rangle = \delta_{\la, \mu}.
\end{eqnarray}
Combining with (\ref{inverse}), we see that the $[\la]$'s with
$\la \in \Pnr$ form a $\C$-linear basis of $\Hn$. So the
restriction of the bilinear form (\ref{form}) to $\Hn \times \Hn$
is a nondegenerate {\em $\C$-valued} bilinear form on $\Hn$ which
again will be denoted by $\langle -, - \rangle$. Let
\begin{eqnarray}  \label{def_H_X}
\Hx =\bigoplus_{n =0}^{\infty} \Hn.
\end{eqnarray}
Then we have an induced non-degenerate bilinear form $\langle -, -
\rangle: \Hx \times \Hx \rightarrow \C$.

\begin{lemma} \label{H1}
\begin{enumerate}
\item [{\rm (i)}] We have the following identifications:
$$rt = \diamondsuit_{0} +\cdots +\diamondsuit_{r-1},\quad
[\Sigma_i] = \diamondsuit_{i-1} - \diamondsuit_{i}, \;1 \le i \le
r-1.$$

\item [{\rm (ii)}] The bilinear form $\langle -, - \rangle$ on
$\mathbb H_1$ is also given by
\begin{eqnarray*}
\langle t, t \rangle = \frac1r, \quad \langle t, [\Sigma_i]
\rangle =0, \quad \langle [\Sigma_i], [\Sigma_j] \rangle = \left
\{
\begin{array}{ll}
       2 &{\rm if} \; i=j,      \\
       -1&{\rm if} \; |i-j| =1, \\
       0 &{\rm otherwise.}
      \end{array}
  \right.
\end{eqnarray*}
\end{enumerate}
\end{lemma}

\begin{proof}
It suffices to verify (i) since (ii) follow from (i) and
(\ref{pairing_Hn}). Applying the localization theorem
(\ref{inverse}) to $\alpha = rt  \in H_T^*(X)$ and $n=1$ yields
\begin{eqnarray*}
 rt = \sum_{i=0}^{r-1} \frac{rt [\xi_i]}{-r^2 t^2} =
\sum_{i=0}^{r-1} \diamondsuit_{i}.
\end{eqnarray*}

Since $\mathfrak I_{i-1}(a_{i-1, 1}, 0) = I_{i}(1:a_{i-1, 1})$
when $1 \le i \le r-1$, the torus $T$ acts on the points $(p_i:
q_i) \in \Sigma_i$ by $s(p_i: q_i) = (p_i: s^{-r}q_i)$. Note that
$\xi_{i-1}$ and $\xi_{i}$ are the points $(1:0)$ and $(0:1)$ in
$\Sigma_i$ respectively. By the localization theorem, we have
\begin{eqnarray*}
1_{\Sigma_i} = \frac{[\xi_{i-1}]}{-rt} + \frac{[\xi_{i}]}{rt} =
-r^{-1}t^{-1} [\xi_{i-1}] + r^{-1}t^{-1} [\xi_{i}] \in
H^*_T(\Sigma_i)'.
\end{eqnarray*}
It follows that $[\Sigma_i] = \diamondsuit_{{i-1}} -
\diamondsuit_{i} \in H^*_T(X)$. This proves (i).
\end{proof}

\subsection{Heisenberg algebra}

In this subsection, we shall generalize the construction of
the Heisenberg algebra in \cite{Vas} (also cf. \cite{Na2}).
Let $i$ be a positive integer, and $Y$ be a
$T$-invariant closed curve in $X$. Define
\begin{align}  \label{Yni}
Y_{n, i} = \{(\xi, \eta) \in X^{[n+i]} \times X^{[n]} \mid \eta
\subset \xi, \,\, \text{Supp}(I_\eta/I_\xi)=\{y\} \in Y \}.
\end{align}
Let $p_1$ and $p_2$ be the projections of $X^{[n+i]} \times
X^{[n]}$ to the two factors respectively. We define the linear
operator $\mathfrak p_{-i}([Y]) \in \text{End}(\Hx)$ by
\begin{eqnarray}   \label{p-i}
\mathfrak p_{-i}([Y])(\alpha) = p_{1}^!(p_2^*\alpha \cup [Y_{n,
i}]) \in \mathbb H_{n+i}
\end{eqnarray}
for $\alpha \in H^{2n}_T(\Xn)$. Note that the restriction of $p_1$
to $Y_{n, i}$ is proper. We define $\mathfrak p_{i}([Y]) \in
\text{End}(\Hx)$ to be the adjoint operator of $\mathfrak
p_{-i}([Y])$. Alternatively, letting $p_2'$ be the projection of
$(X^{[n]})^T \times X^{[n-i]}$ to $X^{[n-i]}$, we see that
\begin{eqnarray*}
\mathfrak p_{i}([Y])(\alpha) = (-1)^{i} \cdot (p_{2}')^! \big (
(\iota_n \times \text{Id})^! \big )^{-1} (p_1^*\alpha \cup
[Y_{n-i, i}]) \in \mathbb H_{n-i}
\end{eqnarray*}
for $\alpha \in H^{2n}_T(X^{[n]})$. Finally, we also put
$\mathfrak p_0([Y])=0$.

Recall
that $\Sigma_0$ is the $z_{0, 2}$-axis in $X_0 \cong \text{\rm
Spec} \, \C[z_{0, 1}, z_{0, 2}]$, and $\Sigma_r$ is the $z_{r-1,
1}$-axis in $X_{r-1} \cong \text{\rm Spec} \, \C[z_{r-1, 1},
z_{r-1, 2}]$. Both $\Sigma_0$ and $\Sigma_r$ are closed in $X$. As
in Lemma~\ref{H1}~(i), using the localization theorem, we obtain
\begin{eqnarray}  \label{0r}
  [\Sigma_0] = -\diamondsuit_{0}, \quad
  \quad  [\Sigma_r] = \diamondsuit_{r-1}.
\end{eqnarray}
Since $\diamondsuit_{0}, \ldots, \diamondsuit_{r-1}$ form a linear
basis of $\mathbb H_1 = H^2_T(X)$, we see from Lemma~\ref{H1}~(i)
and (\ref{0r}) that $\mathbb H_1$ has two more linear bases:
\begin{eqnarray}  \label{two_bases}
\{[\Sigma_0], [\Sigma_1], \ldots, [\Sigma_{r-1}] \}, \qquad
\{[\Sigma_1], [\Sigma_2], \ldots, [\Sigma_r] \}.
\end{eqnarray}

Using one of the two bases in (\ref{two_bases}), we extend by
linearity on $\alpha$ to obtain the linear operator $\mathfrak
p_m(\alpha) \in \End (\Hx)$ for every $\alpha \in \mathbb H_1$.

\begin{theorem} \label{thm_heisenberg}
The linear operators $\mathfrak p_m (\alpha)$, where $m \in \Z$
and $\alpha \in \mathbb H_1 = H^2_T(X)$, satisfy the Heisenberg
algebra commutation relation:
\begin{eqnarray}  \label{thm_heisenberg.0}
[ \mathfrak p_m (\alpha), \mathfrak p_n(\beta) ] = m \delta_{m,-n}
\langle \alpha, \beta \rangle \,\, \text{\rm Id}.
\end{eqnarray}
Furthermore, $\Hx$ is the Fock space (i.e. an irreducible module)
of this Heisenberg algebra with highest weight vector $\vac$ which
denotes the unity in $\mathbb H_0 \subset H^*_T(X^{[0]})$.
\end{theorem}

\begin{proof}
By (\ref{two_bases}), $\{[\Sigma_0], [\Sigma_1], \ldots,
[\Sigma_{r-1}] \}$ is a basis of $\mathbb H_1$. Since $\mathfrak
p_0([\Sigma_i]) = 0$ and $\mathfrak p_{-m}([\Sigma_i])$ is the
adjoint of $\mathfrak p_{m}([\Sigma_i])$, to prove
(\ref{thm_heisenberg.0}) it suffices to prove that
\begin{eqnarray}
 [ \mathfrak p_m ([\Sigma_i]), \mathfrak p_n([\Sigma_j]) ]
  &=& 0   \label{thm_heisenberg.1}  \\
 {[} \mathfrak p_m ([\Sigma_i]), \mathfrak p_{-n}([\Sigma_j]) ]
 &=& m \delta_{m,n} \langle [\Sigma_i], [\Sigma_j] \rangle \,\,
\text{\rm Id}    \label{thm_heisenberg.2}
\end{eqnarray}
for $m, n > 0$ and $0 \le i, j \le r$. When $i \ne j$, $\Sigma_i$
and $\Sigma_j$ are either disjoint or intersect transversally at
exactly one point. Following the argument in \cite{Na2, Vas}, we
conclude that (\ref{thm_heisenberg.1}) and
(\ref{thm_heisenberg.2}) hold when $i \ne j$. To handle the case
$i = j$, we see from Lemma~\ref{H1}~(i) and (\ref{0r}) that
$[\Sigma_i] = - \sum_{0 \le s \le r, s \ne i} [\Sigma_s]$. Thus,
\begin{eqnarray*}
   [\mathfrak p_m ([\Sigma_i]), \mathfrak p_n([\Sigma_i]) ]
 &=& -\sum_{0 \le s \le r, s \ne i} [ \mathfrak p_m ([\Sigma_i]),
   \mathfrak p_n([\Sigma_s]) ] = 0, \\
{[}\mathfrak p_m ([\Sigma_i]), \mathfrak p_{-n}([\Sigma_i])]
 &=& -\sum_{0 \le s \le r, s \ne i} [ \mathfrak p_m ([\Sigma_i]),
\mathfrak p_{-n}([\Sigma_s]) ] \\
&=& m \delta_{m,n} \langle [\Sigma_i], [\Sigma_i] \rangle \,\,
\text{\rm Id}.
\end{eqnarray*}
This completes the proof of (\ref{thm_heisenberg.1}) and
(\ref{thm_heisenberg.2}), and whence (\ref{thm_heisenberg.0}).

To prove the second statement in our theorem, recall that the
classes $[\la]$, as $\la$ runs over all multi-partitions in
$\Pnr$, form a linear basis of $\Hn$. Therefore,
\begin{eqnarray}  \label{eq:eta}
  \sum_{n =0}^{\infty} \dim  (\Hn) \,\, q^n
  = \frac1{\prod_{k=1}^{\infty} (1-q^k)^r}.
\end{eqnarray}
On the other hand, it is well known that the Fock space of the
Heisenberg algebra is irreducible (thanks to the non-degeneracy of
the bilinear form on $\mathbb H_1$), and its character is given by
the right-hand-side of (\ref{eq:eta}). Hence we can identify the
space $\Hx$ with the Fock space of the Heisenberg algebra.
\end{proof}

For latter purpose we make the following definitions. Let $\mu
=(\mu_1, \mu_2, \cdots, \mu_\ell)$ be a partition of the integer
$|\mu |=\mu_1+\cdots+\mu_\ell$, where $\mu_1\geq \cdots \geq
\mu_\ell \geq 1$. We will also make use of another notation $\mu
=(1^{m_1}2^{m_2}\cdots),$ where $m_i$ is the number of parts in
$\mu$ equal to $i$. The length $\ell (\mu)$ is the number $\ell$
or $m_1 +m_2 +\cdots$ in the second notation.
For $0 \le i \le r-1$ and $m\in \Z$, we define
\begin{eqnarray*}
\mathfrak p_{m} (c^i) =\sum_{j=0}^{r-1} \varepsilon^{-ij}
\mathfrak p_m (\diamondsuit_j).
\end{eqnarray*}
For $\la = (\la^0, \cdots, \la^{r-1})$ with $\la^i =
(1^{m_1(i)} 2^{m_2(i)} \ldots)$, define
\begin{eqnarray}
\mathfrak p_{-\la}
 &:=& \prod_{i=0}^{r-1} \mathfrak p_{-\la^i}(\diamondsuit_i) \vac,
 \quad \mathfrak p_{-\la^i}(\diamondsuit_i)
 = \prod_{k \ge 1}
 \mathfrak p_{-k}(\diamondsuit_i)^{m_k(i)}
 \label{p-la}\\
\mathfrak p_{-\la}'
 &:=&  \prod_{i=0}^{r-1}\prod_{k \ge 1}
 \mathfrak p_{-k}(c^i)^{m_k(i)} \vac. \label{p'-la}
\end{eqnarray}

\subsection{Equivariant bundles and equivariant Chern characters}
\label{subsec:chern}

We denote the universal codimension-$2$ subscheme of $\Xn \times
X$ by
$$\mathcal Z_n =\{(\xi,x) \in \Xn \times X \mid x \in
\Supp(\xi) \}. $$
For a line bundle $L$ on $X$, let $L^{[n]}$ denote the
tautological rank-$n$ vector bundle $\pi_{1*}(\mathcal O_{\mathcal
Z_n} \otimes \pi_2^*L)$ on $\Xn$, where $\pi_1$ and $\pi_2$ denote
the projections of $\Xn \times X$ to the factors. When $L$ is
$T$-equivariant over $X$, $L^{[n]}$ is $T$-equivariant over $\Xn$.
This construction is actually valid for every (quasi-)projective
surface (besides $X$).

Next, using \cite{IN} we describe certain distinguished
$T$-equivariant line bundles over $X$, which were earlier defined
in \cite{GSV} by a different and more complicated method. Recall
that $X$ is identified with the subvariety of $(\C^2)^{[r]}$
consisting of the points corresponding to the ideals $I$ in
(\ref{ideals1}) and (\ref{ideals2}) of $\C [z_1, z_2]$. For such
an $I$, $\C [z_1, z_2]/I$ is isomorphic to the regular
representation of $\G$. Denote by $\mathcal O_{\C^2}^{[r]}$ the
tautological rank-$r$ bundle over the Hilbert scheme
$(\C^2)^{[r]}$. Each fiber of the rank-$r$ vector bundle $\big
(\mathcal O_{\C^2}^{[r]} \big )|_X$ over $X$ carries the structure
of the regular representation of $\G$. We obtain $T$-equivariant
line bundles $L_0, L_1, \ldots, L_{r-1}$ over $X$ by decomposing
$\big (\mathcal O_{\C^2}^{[r]} \big )|_X$ according to the
irreducible characters $\g_0, \g_1, \ldots, \g_{r-1}$ of $\G$:
\begin{eqnarray}  \label{def_L_k}
\big (\mathcal O_{\C^2}^{[r]} \big )|_X \cong
\bigoplus_{k=0}^{r-1} \g_k \otimes L_k.
\end{eqnarray}

To understand the fiber $L_k|_{\xi_i}$ of $L_k$ at the $T$-fixed
point $\xi_i$ ($0 \le i \le r-1$), we recall from
(\ref{fixed_pts}) that $\xi_i$ corresponds to the ideal
$(z_1^{i+1}, z_1z_2, z_2^{r-i})$ in $\C [z_1, z_2]$. The fiber of
$(\mathcal O_{\C^2})^{[r]}$ at a point $Z \in ({\C^2})^{[r]}$ is
canonically identified with $H^0(\mathcal O_Z)$. Hence the fiber
of $\big (\mathcal O_{\C^2}^{[r]} \big )|_X$ at $\xi_i$ is
canonically identified as
\begin{eqnarray*}
H^0(\mathcal O_{\xi_i}) = \text{\rm Span } (1, z_1, \ldots,
z_1^{i}, z_2, \ldots, z_2^{r-i-1}).
\end{eqnarray*}
Combining this with the standard action of $\G$ on $\C^2$ and
(\ref{def_L_k}), we conclude that
\begin{eqnarray}  \label{fiber_L_k}
L_k|_{\xi_i} = \left \{
\begin{array}{ll}
  \C \cdot z_2^{r-k} \quad &{\rm if} \,\,\, 0 \le i < k,\\
  \C \cdot z_1^k     \quad &{\rm if} \,\,\, k \le i \le r-1.
\end{array}
  \right.
\end{eqnarray}

Now we study the Chern character of the $T$-equivariant
tautological rank-$n$ vector bundle $L_k^{[n]}$ over $\Xn$. We
begin with the description of its fiber over a $T$-fixed point
$\xi_\la \in \Xn$ where $\la = (\la^0, \ldots, \la^{r-1}) \in
\Pnr$. The fiber of $L_k^{[n]}$ at $\xi_\la$ is
\begin{eqnarray}   \label{fiber_Lkn.0}
L_k^{[n]}|_{\xi_\la} \cong \bigoplus_{i=0}^{r-1}
L_k^{[|\la^i|]}|_{\xi_i^{\la^i}}
\cong \bigoplus_{i=0}^{r-1}\big (L_k|_{\xi_i} \big ) \otimes
H^0(\mathcal O_{\xi_i^{\la^i}}).
\end{eqnarray}
If we write the partition $\la^i$ in terms of its parts as $\la^i
= (\la^i_1, \la^i_2, \cdots, \la^i_\ell)$, then the $\C$-vector
space $H^0(\mathcal O_{\xi_i^{\la^i}}) \subset H^0(\mathcal
O_{X_i})$ has a linear basis
\begin{align*}
\left \{
 1, z_{i,2}, \ldots, z_{i,2}^{\la^i_1-1}, z_{i,1},
z_{i,1} z_{i,2}, \ldots, z_{i,1} z_{i,2}^{\la^i_2-1}, \ldots,
z_{i,1}^{\ell-1}, z_{i,1}^{\ell-1} z_{i,2}, \ldots,
z_{i,1}^{\ell-1} z_{i,2}^{\la^i_\ell-1}
 \right \}
\end{align*}
where $z_{i,1}$ and $z_{i,2}$ are the coordinate functions of
$X_i$. By (\ref{action_on_X_i}), $T$ acts on $z_{i,1}$ and
$z_{i,2}$ by $s(z_{i,1}, z_{i,2}) = (s^rz_{i,1}, s^{-r}z_{i,2})$.
So as a $T$-module,
\begin{eqnarray*}
H^0(\mathcal O_{\xi_i^{\la^i}}) \cong \bigoplus\limits_{\square
\in \la^i} \theta^{rc_{\square}}
\end{eqnarray*}
where $c_{\square}$ is the content of the cell $\square$. By
(\ref{T_action_C^2}), (\ref{fiber_L_k}) and (\ref{fiber_Lkn.0}),
we have
\begin{eqnarray*}
L_k^{[n]}|_{\xi_\la}
 \cong \left ( \bigoplus_{i=0}^{k-1} \, \bigoplus_{\square \in \la^i}
\theta^{(k-r)+rc_{\square}} \right ) \bigoplus \left (
\bigoplus_{i=k}^{r-1} \,\, \bigoplus_{\square \in \la^i}
\theta^{k+rc_{\square}} \right ).
\end{eqnarray*}

Let $\ch_m^T(L_k^{[n]})$ be the $m$-th $T$-equivariant Chern
character of $L_k^{[n]}$. Then,
\begin{eqnarray*}
\ch_m^T(L_k^{[n]})|_{\xi_\la}
 = \frac{1}{m!} \left ( \sum_{i=0}^{k-1}
   \,\, \sum_{\square \in \la^i} \big ( (k-r+rc_{\square})t \big )^m
 + \sum_{i=k}^{r-1} \,\, \sum_{\square \in \la^i} \big (
(k+rc_{\square})t \big )^m  \right ).
\end{eqnarray*}
By the projection formula, we have in $H^*_T(X^{[n]})'$ that
\begin{eqnarray} \label{eq:chern_Lkn}
& &\ch_m^T(L_k^{[n]}) \cup [\xi_\la] = i_\la^!
   \left( \ch_m^T(L_k^{[n]})|_{\xi_\la} \right )
       \\
&=&\frac{1}{m!} \left ( \sum_{i=0}^{k-1} \,  \sum_{\square \in
\la^i}
   \big ( (k-r+rc_{\square})t \big )^m
   + \sum_{i=k}^{r-1} \,\, \sum_{\square \in \la^i} \big (
   (k+rc_{\square})t \big )^m  \right ) [\xi_\la]. \nonumber
\end{eqnarray}

Let $m \ge 0$ and $0 \le k \le r-1$. For a
nonnegative integer $n$, denote
\begin{eqnarray}  \label{ch_km_n}
{\ch}_{k;m}^{[n]} =t^{n-m} \ch_{m}^{T} (L_k^{[n]}) \in \Hn.
\end{eqnarray}
We define a {\em Chern character operator} $\mathfrak G_k $
(respectively, $\mathfrak G_{k;m}$) in $\End (\Hx)$ by sending $a
\in \Hn$ to $a \star \sum_{m \ge 0} {\ch}_{k;m}^{[n]}$
(respectively, to $a \star {\ch}_{k;m}^{[n]}$) in $\Hn$ for each
$n$. Similarly, we define an operator $\mathfrak G_k(z) \in \End
(\Hx)$ by sending $a \in \Hn$ to $a \star \sum_{m \ge 0} z^m
{\ch}_{k;m}^{[n]}$ for each $n$, where $z$ is a variable. By
definition, we have
$$\mathfrak G_k(z) =\sum_{m \ge 0} \mathfrak G_{k;m} z^m.$$

\begin{proposition} \label{prop:chern}
Let $\varsigma (z) =e^{z/2} -e^{-z/2}$. For a given  $\la =(\la^0,
\cdots, \la^{r-1}) \in \Pnr$, we write $\la^i =(\la_1^i, \la_2^i,
\ldots)$ in terms of parts. Then, for $0 \le k \le r-1$, we have
\begin{eqnarray}
\mathfrak G_k(z) \cdot [\la]
 &=& \frac1{\varsigma (rz)}
  \left (
  \sum_{i=0}^{k-1}
  e^{(k-r)z} \left ( \sum_{j=1}^\infty
     e^{(\la^i_j -j+1/2) rz} - \frac1{\varsigma (rz)}  \right )
     \right.   \nonumber \\
 & &   \left.
   \qquad \quad + \sum_{i=k}^{r-1} e^{kz} \left (\sum_{j=1}^\infty
     e^{(\la^i_j -j+1/2) rz} - \frac1{\varsigma (rz)}  \right )
     \right )
\; [\la].  \label{lma_Gzk}
\end{eqnarray}
\end{proposition}
\begin{proof}
Recall that the contents associated to a partition
$\la^i=(\la_1^i, \la_2^i, \ldots, \la^i_\ell)$ are: $1-j, 2-j,
\ldots, \la^i_j -j, 1\le j \le \ell$. A simple algebra
manipulation gives us that
\begin{eqnarray}  \label{eq:alg}
\sum_{\square \in \la^i}
  e^{ c_{\square} z}
  =
  \frac1{\varsigma (z)}
  \left (
   \sum_{j=1}^\infty
     e^{(\la^i_j -j+1/2)z} - \frac1{\varsigma (z)} \right ).
\end{eqnarray}
The formula (\ref{eq:chern_Lkn}) and the definition
(\ref{rescaling}) imply that
\begin{eqnarray*}  \label{eq:basic}
\mathfrak G_k(z) \cdot [\la]
 = \left (
 \sum_{i=0}^{k-1} \, \sum_{\square \in \la^i}
  e^{ (k-r+rc_{\square}) z}
  + \sum_{i=k}^{r-1} \,\, \sum_{\square \in \la^i}
  e^{ (k+rc_{\square})z}
  \right ) \; [\la]
\end{eqnarray*}
which is then easily reduced to (\ref{lma_Gzk}) by using
(\ref{eq:alg}).
\end{proof}
\section{\bf The isomorphisms of rings}
\label{sec:isom}
\subsection{Wreath products and Heisenberg algebra}
\label{sec:wreath}

For a finite group $G$, we denote by $R(G)$ the space of class
functions on $G$, with two distinguished linear bases: one given
by the irreducible characters of $G$ and the other by the
characteristic functions of the conjugacy classes of $G$. We
denote by $R(G)_\Z$ the integral combination of irreducible
characters of $ G$. The standard bilinear form on $R(G )$ is
defined to be such that the irreducible characters form an
orthonormal basis. In this way, $R(G)_\Z$ endowed with this
bilinear form becomes an integral lattice in $R( G )$.

Recall $\G$ is the cyclic group of order $r$. Given a positive
integer $n$, let $\Gamma^n = \Gamma \times \cdots \times \Gamma$
be the $n$-th direct product of $\Gamma$. The symmetric group
$S_n$ acts on $\Gamma^n$ by permutations:
$\sigma (g_1, \cdots, g_n)
  = (g_{\sigma^{ -1} (1)}, \cdots, g_{\sigma^{ -1} (n)}).
$
The wreath product of $\Gamma$ with $S_n$ is defined to be the
semi-direct product
$$
 \Gamma_n = \{(g, \sigma) | g=(g_1, \cdots, g_n)\in {\Gamma}^n,
\sigma\in S_n \}
$$
 with the multiplication
$(g, \sigma)\cdot (h, \tau)=(g \, {\sigma} (h), \sigma \tau ) . $

The conjugacy classes of ${\Gamma}_n$ can be described in the
following way. Let $x=(g, \sigma )\in {\Gamma}_n$, where $g=(g_1,
\cdots, g_n) \in {\Gamma}^n$ and $\sigma \in S_n$. Write
$\sigma $ as a product of disjoint cycles. For each
such cycle $y=(i_1 i_2 \cdots i_k)$ we associate a {\em
cycle-product} $g_{i_k} g_{i_{k -1}} \cdots g_{i_1} \in \Gamma$.
For each integer $i\geq 1$, the number of $k$-cycles in $\sigma$
whose cycle-product equals $a^i$ will be denoted by $m_k(i)$.
Denote by $\la^i$ the partition $(1^{m_1 (i)} 2^{m_2 (i)} \cdots
)$. Then each element $x=(g, \sigma)\in {\Gamma}_n$ gives rise to
a multi-partition
\begin{eqnarray*}
( \la^0,\cdots, \la^{r-1}) \in \Pnr
\end{eqnarray*}
which will be referred to as the {\em type} of $x$. It is well
known (cf. \cite{Mac}) that two elements of ${\Gamma}_n$ are
conjugate in ${\Gamma}_n$ if and only if they have the same type.

Given a partition $\mu = (1^{m_1} 2^{m_2} \ldots )$, we define
$z_{\mu} = \prod_{i\geq 1}i^{m_i}m_i!$. The order of the
centralizer of an element $x = (g, \sigma) \in {\Gamma}_n$ of
type $\la =( \la^0,\cdots, \la^{r-1})$ is
\begin{eqnarray*}
Z_{\la}= r^{\ell(\la)} \prod\limits_{i=0}^{r-1} z_{\la^i}.
\end{eqnarray*}

For $\mu \vdash m$, we associate with the irreducible
character $s_\mu$ of the symmetric group $S_m$ and its
corresponding representation $U_\mu$. Denote by $V_i$ the
irreducible representation of $\G$ associated to the character
$\g_i$. Then the wreath product $\G_m =\G^m \rtimes S_m$ acts
irreducibly on $U_\mu \otimes V_i^{\otimes m}$ where $\G^m$ acts
on the second tensor factor only and $S_m$ acts diagonally. We
denote the associated character of $\G_m$ by $s_\mu(\g_i)$.

More generally, let $\la =( \la^0,\cdots, \la^{r-1})$ be a
multi-partition with $\| \la\| := |\la^0| + \ldots + |\la^{r-1}|
=n$. Naturally $W_\la = \bigotimes_i \left(U_{\la^i} \otimes
V_i^{\otimes |\la^i|}\right)$ is a representation of
$\prod_i \G_{|\la^i|}$ which is a subgroup of $\Gn$.
The induced $\Gn$-representation
$\text{Ind}^{\Gn}_{\prod_i \G_{|\la^i|} } W_\la$
is irreducible and its associated character is given by
\begin{eqnarray}   \label{eq:product}
 s_\la = \prod_i s_{\la^i}(\g_i).
\end{eqnarray}
Letting $d_{\la^i}$ be the degree of $s_{\la^i}$ (i.e. the
dimension of the corresponding representation) and
$h(\la) =\prod_i h(\la^i)$ where $h(\la^i)$ is the hook
product associated to $\la^i$, we see that
the degree $d_\la$ of the character $s_\la$ is
\begin{eqnarray}  \label{eq:hook}
d_\la = \frac{n!}{\prod_i |\la^i| !} \prod_i d_{\la^i}
  = \frac{n!}{\prod_i h(\la^i)}
  = \frac{|\Gn|}{r^n  h(\la)}.
\end{eqnarray}

We define
\[
  \RG = \bigoplus_{n= 0}^\infty R({\Gamma}_n).
\]
It carries a bilinear form, denoted by $\langle -,- \rangle$,
induced from the standard one on each $R({\Gamma}_n)$. There is an
action on $\RG$ (cf. \cite{FJW, Wa1}) of a Heisenberg algebra
which is generated by $\mathfrak a_n (\g_i), n \in \Z, 0 \le i \le
r-1$ with the commutation relation:
\begin{eqnarray} \label{eq:orth}
[ \mathfrak a_m (\g_i), \mathfrak a_n (\g_j)]
 = m \delta_{m,-n} \delta_{i,j} \text{Id}.
 \end{eqnarray}
Here $\text{Id}$ denotes the identity operator on $\RG$. We
further extend the definition of $\mathfrak a_n (\g)$ to all $\g
\in R(\G)$ by linearity. The operators $\mathfrak a_{ n} (\g)$
(respectively, $\mathfrak a_{ -n} (\g)$) for $n>0$ have been
defined in terms of restriction functors (respectively, induction
functors) and will be referred to as annihilation operators
(respectively, creation operators). It is known (cf. \cite{FJW})
by (\ref{eq:product}) and the formulation of the Heisenberg
operators that, for $m>0$ and $\la =(\la^0, \cdots, \la^{r-1})$,
\begin{eqnarray}
 \mathfrak a_{\pm m} (\g_j) ( s_{\la^i}(\g_i) )
 &=& 0, \quad \quad\quad j \neq i  \label{eq:simple1} \\
 \mathfrak a_{ - m} (\g_i) ( s_{\la^i}(\g_i) )
  &=& \sum_{|\mu| = m +|\la^i|} c(m, \la^i; \mu) \; s_\mu (\g_i)
   \label{eq:simple2}
\end{eqnarray}
where $c(m, \la^i; \mu) \in \mathbb Q$ depends only on $m$ and the
partitions $ \la^i, \mu$. It follows that
\begin{eqnarray}
 \mathfrak a_{ - m} (\g_i) ( s_{\la})
  &=& \sum_{|\mu| = m +|\la^i|} c(m, \la^i; \mu) \;
  s_{(\la^0, \cdots, \la^{i-1}, \mu, \la^{i+1},\cdots,
  \la^{r-1})}.
 \label{eq:simple3}
\end{eqnarray}
In particular, the coefficients $c(m, \la^i; \mu)$'s do not depend
on $r$, i.e. they are the same structure constants as for the
symmetric group case.

Denote by $c^i$ the class function in $R(\G)$ which takes value
$r$ on the element $a^i$ and $0$ elsewhere. Following \cite{FJW},
we can use the character table of $\G$ and linearity to define
another set of generators $\mathfrak a_m (c^i), 0 \le i \le r-1$,
of the Heisenberg algebra by
 $$
 \mathfrak a_m (c^i)
  =\sum_{k=0}^{r-1} \varepsilon^{-ik} \mathfrak a_m (\g_k).
 $$
We have
$$
[ \mathfrak a_m (c^i), \mathfrak a_n (c^j)]
 = rm \,  \delta_{m,-n} \delta_{i+j,r} \; \text{Id}.$$
Now the space $\RG$ is an irreducible module over this Heisenberg
algebra. Denote by $\vac$ the unity in $R(\G_0)$. Then we have
$$\RG = \C [ \mathfrak a_{-n} (c^i) \mid n>0, 0 \le i \le
r-1 ] \cdot \vac.$$

Given $\la =( \la^0,\ldots, \la^{r-1})$ where $\la^i =(1^{m_1(i)}
2^{m_2(i)} \ldots)$ and $ \| \la\| =n$, we define
\begin{eqnarray*}
\mathfrak a_{-\la} &=& \prod_{i=0}^{r-1} \prod_{k \ge 1}
   \mathfrak a_{-k}(c^i)^{m_k(i)} \vac;    \\
{\mathfrak a}^R_{-\la} &=& \prod_{i=0}^{r-1} \prod_{k \ge 1}
   \mathfrak a_{-k}(\g_i)^{m_k(i)} \vac.
\end{eqnarray*}

\begin{remark}  \label{rem:cc}
We can show that $\mathfrak a_{-\la}$ is the class function on
$\Gn$ which takes value $Z_\la$ on the conjugacy class of type
$\la$, and $0$ elsewhere.
\end{remark}

We have three distinguished linear bases for $R(\Gn)$:
\begin{eqnarray*}
\qquad\qquad  \{\mathfrak a_{-\la} \}, \{\mathfrak a^R_{-\la}\},
\{s_\la\}, \qquad \la \in \Pnr.
\end{eqnarray*}
The transition matrix between
the bases $\{\mathfrak a^R_{-\la}\}$ and $\{s_\la\}$ is given by
the $r$-tuple version of the character table matrix of the
symmetric groups. The transition matrix between the bases
$\{\mathfrak a_{-\la}\}$ and $\{s_\la\}$ is precisely the
character table of $\Gn$. 

\subsection{The class algebras of wreath products}
\label{sec:modify}

For any given finite group $G$, the space $R(G)$ of complex-valued
class function carries an algebra structure (which is often
refereed to as the {\em class algebra} of $G$) given by the {\em
convolution product}:
\[
  (\beta \circ \g) (x) =  \sum_{y \in G} \beta(x y^{-1}) \g(y),
  \quad \beta, \g \in R(G), \; x \in G.
\]
It is well known that, for two irreducible characters $\beta$ and
$\gamma$ of $G$,
\begin{eqnarray}  \label{eq:idem}
   \beta  \circ \g=  \delta_{\beta,\g } \frac{|G|}{d_\g} \; \g,
\end{eqnarray}
where $d_{\g}$ is the {\em degree} of the irreducible character
$\g$. In this paper, we mainly apply these considerations to the
group $\Gn$.

Let $x$ be an element of $\Gn$ of type $\la =(\la^0, \cdots,
\la^{r-1}) \in \Pnr$. We define the {\em modified type} of $x$ to
be $\widetilde{\la} \in \mathcal P_{n- \ell}(r)$, where $\ell
=\ell(\la^0)$, as follows: $\widetilde{\la}^i =\la^i$ for $i > 0$
and $\widetilde{\la}^0 = (\la^0_1 -1, \cdots, \la^0_\ell -1)$ if
we write the partition ${\la}^0= (\la^0_1, \cdots, \la^0_\ell)$.
The {\em degree} $\|x\|$ of $x$ is defined to be  $\|
\widetilde{\la}\|$, which equals $\| \la \| - \ell (\la^0)$.
Apparently, two conjugate elements have the same degree, and it
makes sense to talk about the degree of a conjugacy class of
$\Gn$. It was shown \cite{Wa3} that
$$\|x y\|  \leq \|x\| +\|y\|, \qquad x, y \in \Gn  $$
and thus the degree defines a ring filtration for $R(\Gn)$. Also
it was shown that this notion of degree coincides with the
fermionic degree number or age introduced by Zaslow and Ito-Reid.
Denote the associated graded ring by
$$\mathcal G_\G (n) =\bigoplus_{i=1}^n
\mathcal G_\G^i(n)$$
where $\mathcal G_\G^i(n)$ consists of elements in $R(\Gn)$ of
degree $i$.

Among the conjugacy classes of $\Gn$, the conjugacy class
$K_{[21^{n-2}]}$ of ``transposition" plays a special role. We
decompose $K_{[21^{n-2}]}$ into a disjoint union:

\begin{eqnarray*}
  K_{[21^{n-2}]}
  = \bigsqcup_{j =1}^n K_{[21^{n-2}]}(j),
  \qquad K_{[21^{n-2}]}(j)
  = \bigsqcup_{1 \le i<j} K_{[21^{n-2}]}(i,j)
\end{eqnarray*}
where $K_{[21^{n-2}]}(i,j)$ consists of elements $( (g_1, \ldots,
g_n), (i,j)) \in \Gn$, where all $g_k$ except $g_i$ and $g_j$ are
equal to $1\in \G$, and $g_j = g_i^{-1}$  runs over $\G$.

The elements $M_j =M_{j;n} = \sum_{x \in K_{[21^{n-2}]}(j)} x,$
\,$ 1 \le j \le n,$  in the group algebra $\C[\Gn]$ were
introduced in \cite{Wa2}, and they satisfy several favorable
properties, e.g. $M_i M_j =M_j M_i$ for all $i,j$. Note that
$M_1=0$. We call them the Jucys-Murphy (JM) elements for $\Gn$
since they reduce to the usual Jucys-Murphy elements of the
symmetric group $S_n$ when $\G$ is trivial.

Let $\G^{(i)}$ denote the $i$-th copy of $\G$ in $\Gn$. Given
$\alpha \in R(\G)$, we denote $\alpha^{(i)}$ to be the copy of
$\alpha$ in $\C[\G^{(i)}] \subset \C[\Gn]$.
Letting $z$ be a formal parameter, we define
\begin{eqnarray*}
\Xi_n^m(\alpha) = \frac1{r^m m!} \sum_{i=1}^n  {M_i}^m \circ
\alpha^{(i)},
\qquad \Xi_{n}^{(\alpha)}(z) = \sum_{m\ge 0} \Xi_n^m(\alpha)\;
{z^m}.
\end{eqnarray*}
A key property of
$\Xi_n^m(\alpha)$ is that $\Xi_n^m (\alpha)$ lies in $R(\Gn)$. We
further define the operator $\mathfrak O^m(\alpha) \in \End (\RG)$
(respectively, $\mathfrak O^{(\alpha)}(z )$) to be the convolution
product with $\Xi_n^m(\alpha)$ (respectively,
 $\Xi_{n}^{(\alpha)}(z)$) in $R(\Gn)$ for every $n
\ge 0$.

Let us introduce the following vertex operator associated to
$\g_i, 0\le i \le r-1$:
\begin{eqnarray*}
V(\g_i;w,q)
 &=& \exp \left(\sum_{k > 0} \frac{(q^{k}-1)w^{k}}{k} \mathfrak
a_{-k}(\g_i)
          \right)
 \exp \left(\sum_{k > 0} \frac{(1-q^{-k})w^{-k}}{k} \mathfrak a_k(\g_i)
      \right).
\end{eqnarray*}
Write $V(\g_i;w,q) = \sum_{m \in \Z} V_m(\g_i;q)w^{-m}.$ It is
established in \cite{Wa2}\footnote{Our convention on $\mathfrak
O^{(\alpha)}(z )$ here differs from \cite{Wa2} by a scaling of $z$
by a factor $r$.} that
\begin{eqnarray}  \label{eq:vo}
 \mathfrak O^{(\g_i)}(z )
 =\frac{e^{r z}}{(e^{r z}-1)^2} \left( V_0(\g_i; e^{r z})-1 \right),
 \quad 0 \le i \le r-1.
\end{eqnarray}
When $\G$ is trivial, this specializes to a result of \cite{LT}.

\begin{proposition}  \label{prop:JM}
For $\la =(\la^1, \cdots, \la^{r-1}) \in \Pnr$ with $\la^i
=(\la^i_1, \la^i_2, \ldots)$, we have
\begin{eqnarray} \label{eq:multi}
\mathfrak O^{(\g_i)}(z )  \cdot s_\la
 =\frac1{\varsigma (r z)}
 \left(\sum_{j=1}^\infty e^{ (\la^i_j -j+1/2)r z}
  -\frac1{\varsigma (r z)}\right) \; s_\la.
\end{eqnarray}
\end{proposition}

\begin{proof}
First let $r=1$. It is known that the eigenvalues of the JM
elements of the symmetric group $S_n$ on an irreducible module
$V_\nu$ (where $\nu$ is a partition of $n$) are given by the
contents of the Young tableaux associated to $\nu$, cf. \cite{LT,
LQW2}. Then a simple calculation using generating function gives
us (\ref{eq:multi}). By (\ref{eq:vo}),
\begin{eqnarray} \label{eq:special}
\frac{e^{z}\left( V_0(1; e^{z})-1 \right) }{(e^{z}-1)^2}  \cdot
s_\nu
 =\frac1{\varsigma (z)}
 \left(\sum_{j=1}^\infty e^{(\nu_j -j+1/2) z}
  -\frac1{\varsigma (z)}\right) \; s_\nu.
\end{eqnarray}

Recall that $s_\la = \prod_j s_{\la^j} (\g_j)$. So by applying
(\ref{eq:simple1}), (\ref{eq:vo}) and (\ref{eq:special}) (with $z$
replaced by $r z$, $\nu$ by $\la^i$, etc) to the $\g_i$-component,
we obtain (\ref{eq:multi}).
\end{proof}

\subsection{An isomorphism in the equivariant setup}

Recall the well-known fact that the Hilbert-Chow morphism $\tau_n:
\Xn \rightarrow X^n/S_n$, which sends an element of $\Xn$ to its
support, is a resolution of singularities. Combining with the
minimal resolution $\pi: X \rightarrow \C^2/\G$, we have the
following commutative diagram
$$\CD X^{[n]} @>>> X^n  /S_n \\ @VVV @VVV \\ \C^{2n}  /\Gn @<{\cong}<<
(\C^2/\G)^n /S_n \endCD $$
which give us a (crepant) resolution of singularities $\pi_n: \Xn
\rightarrow \C^{2n}  /\Gn$. This was the starting point on why
Hilbert schemes have something to do with wreath products (cf.
\cite{Wa1} and the references therein).

We define a linear map $\phi: \Hn \longrightarrow R(\Gn)$ by
letting
\begin{eqnarray} \label{def:phi}
\phi ([\la]) = s_\la
\end{eqnarray}
for each $\la \in \Pnr$, which is clearly an isomorphism of vector
spaces. Putting these isomorphisms together, we obtain a linear
isomorphism $\phi: \Hx \to \RG$.

\begin{lemma}  \label{lma:phi_p}
The linear isomorphism $\phi: \Hx \to \RG$ commutes with the
action of the Heisenberg creation operators. More precisely, for
$m>0$, we have
\begin{eqnarray*}
\phi \big (\mathfrak p_{-m}([\Sigma_i])\cdot [\la] \big )
 = \left \{
\begin{array}{ll}
  -\mathfrak a_{-m}(\g_0) \cdot s_\la  \label{phi_p1}
      \quad &{\rm if} \,\,\, i = 0,   \\
  \mathfrak a_{-m}(\g_{i-1} -\g_i) \cdot s_\la
      \quad &{\rm if} \,\,\, 1 \le i \le r-1,  \\
  \mathfrak a_{-m}(\g_{r-1}) \cdot s_\la
      \quad &{\rm if} \,\,\, i=r.
\end{array}
  \right.
\end{eqnarray*}
\end{lemma}

\begin{proof}
When $r = 1$, $T$ acts on $\C^2$ by $s(z_1, z_2) = (sz_1,
s^{-1}z_2)$. In this case, Lemma~\ref{lma:phi_p} follows from the
results in \cite{Vas}, i.e., given a partition $\la$, we have
\begin{eqnarray}
   -\mathfrak p_{-m}([\Sigma_0])[\la]
&=&\sum_{|\mu| = m +|\la|} c(m, \la; \mu)\, [\mu],
\label{lma:phi_p.1}
\end{eqnarray}
where the coefficients $c(m, \la; \mu)$ were defined in
(\ref{eq:simple2}). We note that (\ref{lma:phi_p.1}) remains to be
valid if the $T$-action on $\C^2$ is replaced by $s(z_1, z_2) =
(s^kz_1, s^{-k}z_2)$ with $k> 0$ and $[\la]$ is defined to be
$h(\la)^{-1} (-kt)^{-|\la|} [\xi_{\la}]$ as in (\ref{rescaling}).

Assume $r = 2$. Let $\la=(\la^0, \la^{1}) \in {\mathcal P}_n(2)$.
Recall that $\Sigma_0 \subset X_0$ and that $\xi_0$ is the only
$T$-fixed point on $\Sigma_0$ and $X_0$. Now the inclusion map
$\iota_0: X_0 \to X$ is $T$-equivariant. It induces the
$T$-equivariant inclusion map $\iota_{01}: X_0^{[\tilde n+m]} \to
X^{[n+m]}$ sending $\xi$ to $\xi + \xi_1^{\la^1}$, where $\tilde n
= |\la^0|$. Similarly, we have $T$-equivariant morphisms
$\iota_{02}: X^{[n+m]} \times X_0^{[\tilde n]} \to X^{[n+m]}
\times \Xn$ and $\iota_{03}: X_0^{[\tilde n +m]} \times
X_0^{[\tilde n]} \to X^{[n+m]} \times X_0^{[\tilde n]}$ by adding
$\xi_1^{\la^1}$ suitably. By (\ref{p-i}) and the projection
formula, we obtain
\begin{eqnarray*}
   \mathfrak p_{-m}([\Sigma_0])[\xi_{\la}]
&=&p_{1}^! \big (p_2^*[\xi_{\la}] \cup [(\Sigma_0)_{n, m}] \big )
   = p_{1}^! \big ([X^{[n+m]} \times \xi_{\la}]
   \cup [(\Sigma_0)_{n, m}] \big )   \\
&=&p_{1}^! \big ((\iota_{02})^![X^{[n+m]} \times \xi_0^{{\la}^0}]
   \cup [(\Sigma_0)_{n, m}] \big ) \\
&=&p_{1}^!({\iota}_{02})^! \big ([X^{[n+m]} \times
\xi_0^{{\la}^0}]
   \cup ({\iota}_{02})^*[(\Sigma_0)_{n, m}] \big )   \\
&=&p_{1}^!({\iota}_{02})^! \big ([X^{[n+m]} \times
\xi_0^{{\la}^0}]
   \cup ({\iota}_{03})^![(\Sigma_0)_{\tilde n, m}] \big )  \\
&=&p_{1}^!({\iota}_{02})^!({\iota}_{03})^! \big (
   ({\iota}_{03})^*[X^{[n+m]} \times \xi_0^{{\la}^0}]
   \cup [(\Sigma_0)_{\tilde n, m}] \big )   \\
&=&p_{1}^!({\iota}_{02})^!({\iota}_{03})^! \big (
   [X_0^{[\tilde n+m]} \times \xi_0^{{\la}^0}]
   \cup [(\Sigma_0)_{\tilde n, m}] \big )
\end{eqnarray*}
where $p_1$ and $p_2$ are the projections of $X^{[n+m]} \times
X^{[n]}$ to the two factors, and $(\Sigma_i)_{\tilde n, m} \subset
X_0^{[\tilde n +m]} \times X_0^{[\tilde n]}$ is defined similarly
as in (\ref{Yni}). Let ${\tilde p}_1$ and ${\tilde p}_2$ be the
projections of $X_0^{[\tilde n+m]} \times X_0^{[\tilde n]}$ to the
two factors. Then
\begin{eqnarray}  \label{reduction}
 \mathfrak p_{-m}([\Sigma_0])[\xi_{\la}]
   &= & p_{1}^!({\iota}_{02})^!({\iota}_{03})^! \big (
  ({\tilde p}_2)^*[\xi_{\la^0}] \cup [(\Sigma_0)_{\tilde n, m}]
  \big )  \\
 &=&({\iota}_{01})^! ({\tilde p}_{1})^!\big (
  ({\tilde p}_2)^*[\xi_{\la^0}] \cup [(\Sigma_0)_{\tilde n, m}]
  \big )  \nonumber \\
 &=& ({\iota}_{01})^! \,\, \mathfrak
 p_{-m}([\Sigma_0])[\xi_{\la^0}].  \nonumber
\end{eqnarray}
Combining this with (\ref{lma:phi_p.1}), we conclude that
\begin{eqnarray}   \label{r=2;i=0}
 \mathfrak p_{-m}([\Sigma_0])[\la]
 = -\sum_{|\mu| = m +|\la^0|} c(m, \la^0; \mu) [(\mu, \la^1)].
\end{eqnarray}
Similarly, we can show that
\begin{eqnarray}  \label{r=2;i=2}
 \mathfrak p_{-m}([\Sigma_2])[\la]
 = \sum_{|\mu| = m +|\la^1|} c(m, \la^1; \mu) [(\la^0, \mu)].
\end{eqnarray}
Recall that $[\Sigma_0] + [\Sigma_1] + [\Sigma_2] = 0$ since $r =
2$. So we have
\begin{align}  \label{r=2;i=1}
 \mathfrak p_{-m}([\Sigma_1])[\la]=
 \sum_{|\mu| = m +|\la^0|} c(m, \la^0; \mu) [(\mu, \la^1)] -
\sum_{|\mu| = m +|\la^1|} c(m, \la^1; \mu) [(\la^0, \mu)].
\end{align}
Thus Lemma~\ref{lma:phi_p} holds for $r = 2$ by comparing with
(\ref{eq:simple3}). Again, we stress that (\ref{r=2;i=0}),
(\ref{r=2;i=2}) and (\ref{r=2;i=1}) are still valid if the
$T$-action on $X$ is modified so that $T$ acts on $X_i$ ($i =
0,1$) by $s(z_{i,1}, z_{i,2}) = (s^kz_{i,1}, s^{-k}z_{i,2})$ for a
$k > 0$.

Now assume $r > 2$. Let $\la=(\la^0, \ldots, \la^{r-1})
\in \Pnr$. The same arguments as in the proof of (\ref{r=2;i=0})
show that Lemma~\ref{lma:phi_p} is true for $i = 0$ or $r$.
Let $1 \le i \le (r-1)$. Then the $T$-fixed points in the
projective line $\Sigma_i$ are $\xi_{i-1}$ and $\xi_i$. Let
$\widetilde{X} = X_{i-1} \cup X_i$ be equipped with the induced
$T$-action. Then, $\xi_{i-1}$ and $\xi_i$ are the only $T$-fixed
points in $\widetilde{X}$, and we can apply the results in the
previous paragraph to $\widetilde{X}$. The $T$-equivariant
inclusion map $\tilde \iota: \widetilde{X} \to X$ induces the
$T$-equivariant inclusion map ${\tilde \iota}_1:
\widetilde{X}^{[\tilde n+m]} \to X^{[n+m]}$ which sends $\xi$ to
$\xi + \sum_{0 \le j \le r-1, \, j \ne i-1, i} \xi_j^{\la^j}$,
where $\tilde n = |\la^{i-1}|+|\la^{i}|$. Put $\tilde \la =
({\la}^{i-1}, {\la}^i) \in {\mathcal P}_{\tilde n}(2)$ and
$\xi_{\tilde \la} = \xi_{i-1}^{\la^{i-1}} + \xi_{i}^{\la^{i}} \in
\widetilde{X}^{[\tilde n]}$. Arguments parallel to the proof of
(\ref{reduction}) show that
\begin{eqnarray*}
\mathfrak p_{-m}([\Sigma_i])[\xi_{\la}] = ({\tilde \iota}_1)^!
\,\, \mathfrak p_{-m}([\Sigma_i])[\xi_{\tilde \la}].
\end{eqnarray*}
Combining this with (\ref{r=2;i=1}), we conclude that
\begin{eqnarray*}
  \mathfrak p_{-m}([\Sigma_i])[\la]
 &=&\sum_{|\mu| = m +|\la^{i-1}|}
c(m, \la^{i-1}; \mu) [(\la^0, \ldots,
   \la^{i-2}, \mu, \la^i, \la^{i+1}, \ldots, \la^{r-1})] \nonumber \\
 && - \sum_{|\mu| = m +|\la^i|} c(m, \la^i; \mu) [(\la^0, \ldots,
\la^{i-2},
   \la^{i-1}, \mu, \la^{i+1}, \ldots, \la^{r-1})].
\end{eqnarray*}
Now we complete the proof of Lemma~\ref{lma:phi_p} for $r > 2$ by
comparing with (\ref{eq:simple3}).
\end{proof}

\begin{theorem}  \label{th:isom1}
There is a canonical ring isomorphism $\phi: \Hn \longrightarrow
R(\Gn)$ satisfying
\begin{eqnarray}
\phi ([\la]) &=& s_\la  \nonumber
\\
\phi (\mathfrak p_{-\la}) &=& \mathfrak a^R_{-\la}
     \label{eq:precise2} \\
\phi (\mathfrak p_{-\la}') &=& \mathfrak a_{-\la}
     \label{eq:precise3}
\end{eqnarray}
for each $\la \in \Pnr$. In addition, $\phi$ is an isometry with
respect to the standard bilinear forms on $\Hx$ and $\RG$, and
commutes with the action of Heisenberg algebra.
\end{theorem}
\begin{proof}
We have noted that the map $\phi$ defined by (\ref{def:phi}) is a
linear isomorphism. The map $\phi$ is an isometry since for two
arbitrary $\la$ and $\mu$ in $\Pnr$, we have
\begin{eqnarray*}
\langle [\la], [\mu] \rangle = \delta_{\la,\mu} = \langle s_\la,
s_\mu \rangle.
\end{eqnarray*}
The compatibility of $\phi$ with the Heisenberg
creation operators was verified in Lemma~\ref{lma:phi_p}. Since the
annihilation operators are the adjoints of the creation
operators with respect to the bilinear forms, they are also
compatible with $\phi$.

By (\ref{xi_cup}), (\ref{rescaling}) and the definition of
$\star$, we obtain
\begin{eqnarray*}
[\la] \star [\mu] = \delta_{\la,\mu} r^n h(\la) \; [\la].
\end{eqnarray*}
On the other hand, it follows from (\ref{eq:hook}) and
(\ref{eq:idem}) that
\begin{eqnarray*}
s_\la \circ s_\mu = \delta_{\la,\mu} r^n h(\la) \; s_\la.
\end{eqnarray*}
Thus $\phi: \Hn \longrightarrow R(\Gn)$ is actually a ring
isomorphism.

Lemma~\ref{lma:phi_p} implies that for each $m>0$ and $\la \in
\Pnr$, we have
\begin{eqnarray}  \label{eq:clear}
\phi \big (\mathfrak p_{-m}(\diamondsuit_i )\cdot [\la] \big )
 =  \mathfrak a_{-m}(\g_i) \cdot s_\la, \qquad  0\le i \le r-1.
\end{eqnarray}
Now (\ref{eq:precise2}) follows from an induction argument
and the definitions of $\mathfrak p_{-\la}$
and $\mathfrak a^R_{-\la}$. Finally, note from definitions that
the transition matrix between
$\mathfrak p_{-\la}$ and $\mathfrak p_{-\la}'$ coincides with the
one between $\mathfrak a^R_{-\la}$ and $\mathfrak a_{-\la}$. Thus
(\ref{eq:precise3}) follows from (\ref{eq:precise2}).
\end{proof}

\subsection{An isomorphism for ordinary cohomology rings}

By Lemma~\ref{H1}~(i), we have $rt = \sum_{i=0}^{r-1} \diamondsuit_i$.
On the other hand, $c^0 =\sum_{i=0}^{r-1} \g_i$. It follows from
(\ref{eq:clear}) that, for each $m>0$ and $\la \in \Pnr$,
\begin{eqnarray} \label{eq:unit}
\phi \big (\mathfrak p_{-m}(rt)\cdot [\la] \big )
 =  \mathfrak a_{-m}(c^0) \cdot s_\la.
\end{eqnarray}

Denoting $\la =(\la^0, \cdots, \la^{r-1})$ and $\la^i =
(1^{m_1(i)}2^{m_2(i)} \cdots)$, we define
\begin{eqnarray} \label{b_la}
\mathfrak b_{-\la}
 &= & \prod_{k \ge 1} \left( \mathfrak a_{-k} (c^0)^{m_k(0)}
 \prod_{i=1}^{r-1} \mathfrak a_{-k}
(\g_{i-1}-\g_i)^{m_k(i)} \right) \vac \in \RG \\
\mathfrak q_{-\la}^T
 &= & \prod_{k \ge 1}
 \left ( \mathfrak p_{-k} (rt )^{m_k(0)}
  \prod_{i=1}^{r-1} \mathfrak p_{-k} ([\Sigma_i])^{m_k(i)}
  \right)
  \vac \in \Hx. \nonumber
\end{eqnarray}
It follows by induction from Lemma~\ref{lma:phi_p} and
(\ref{eq:unit}) that
\begin{eqnarray}   \label{eq:graded}
 \phi (\mathfrak q_{-\la}^T )
 =  \mathfrak b_{-\la}.
\end{eqnarray}

Note that a linear basis for the ordinary cohomology $H^*(X)$ is
given by
\begin{eqnarray*}
1_X \in H^0(X), \qquad \Sigma_1, \ldots, \Sigma_{r-1} \in
H^2(X).
\end{eqnarray*}
Using the usual Heisenberg algebra construction
\cite{Na1} (in Nakajima's notations, the Heisenberg operators were
denoted by $P_{1_X}[m]$ and $P_{\Sigma_i}[m]$), we can construct a
linear basis $\{ \mathfrak q_{ -\la} \}_{\la \in \Pnr}$ for
$H^*(\Xn)$, where we have denoted
\begin{eqnarray}   \label{q_la}
\mathfrak q_{-\la} =
 \prod_{k \ge 1}
 \left( P_{r \cdot 1_X}[-k]^{m_k(0)} \prod_{i=1}^{r-1}
 P_{\Sigma_i}[-k]^{m_k(i)}
 \right) \vac
\end{eqnarray}
if $\la =(\la^0, \ldots, \la^{r-1})$ with $\la^i =
(1^{m_1(i)}2^{m_2(i)} \cdots)$. Using the degrees of the
Heisenberg operators (cf. \cite{Na2}), the cohomology degree of
$\mathfrak q_{-\la}$ is computed to be
$$\sum_{k \ge 1} \left (
 2(k-1) m_k(0) +\sum_{i=1}^{r-1} (2(k-1)+2) m_k(i) \right )
 =2 \big(\| \la \| - \ell(\la^0) \big).$$

Let $F^p = \sum_{k \le p} t^{n-k} H^{2k}_T (\Xn)$. Then
\begin{eqnarray*}
F^0 \subset F^1 \subset \ldots \subset F^n =\Hn
\end{eqnarray*}
defines a
filtration on the ring $\Hn$. By a spectral sequence argument (cf.
\cite{Na3, Vas}), $H^{*}(\Xn)$ is identified with the graded ring
of $\Hn$ associated to this filtration, and $\mathfrak q_{-\la}
\in H^{*}(\Xn)$ is the element associated to $\mathfrak q_{-\la}^T
\in \Hn$.

On the other hand, note that $(\g_j -\g_{j+1})$ is homogeneous in
$\mathcal G_\G(1)$ of degree $1$ and $c^0$ is of degree 0. By the
definition of the degree of a conjugacy class of $\Gn$ in
Sect.~\ref{sec:modify} and Remark~\ref{rem:cc}, the degree of
$\mathfrak b_{-\la}$ equals $\| \la \| - \ell(\la^0).$

Thus, by (\ref{eq:graded}) above, the filtration on the ring $\Hn$
is compatible (up to a multiple $2$) with the degree filtration on
$R(\Gn)$ defined in Sect.~\ref{sec:modify}. We can regard
$\mathfrak b_{-\la}$ as an element in $\mathcal G^{*}_\G (n)$.
Further note that $c^0$ is $r$ times the identity element of the
ring $\mathcal G^{*}_\G (1)$, which is compatible with $r\cdot
1^X$ in $H^0(X)$. Also compare (\ref{b_la}) and (\ref{q_la}).
Summarizing, we have established the following.

\begin{theorem}  \label{th:isom2}
Rescale the grading on $\mathcal G^{*}_\G (n)$ by a multiple of
$2$. Then, the map from $H^{*}(\Xn)$ to $\mathcal G^{*}_\G (n)$
which sends $\mathfrak q_{ -\la}$ to $\mathfrak b_{-\la}$ for each
$\la \in \Pnr$ is a graded ring isomorphism. $\qed$
\end{theorem}

\begin{remark}
When $r=1$, i.e. $\G$ is trivial, $X$ is the affine plane and
$\Gn$ is the symmetric group $S_n$. In this case, our
Theorem~\ref{th:isom1} specializes to a theorem in \cite{Vas}, and
Theorem~\ref{th:isom2} specializes to a theorem in \cite{LS, Vas}
(also cf. \cite{LQW1}). Theorem~\ref{th:isom2} further supports
and enhances Ruan's conjecture on the existence of an isomorphism
between the cohomology ring of a hyperkahler resolution and
Chen-Ruan orbifold cohomology ring, \cite{Ru} (also cf.
\cite{Wa3}). Ruan's conjecture does not provide any explicit map
which realizes a ring isomorphism. As a special case of a more
general theorem in \cite{EG}, the ring isomorphism in
Theorem~\ref{th:isom2} without the explicit map has been verified
in a very indirect way.
\end{remark}

\begin{remark}
The structure constants of the algebra $\mathcal G_\G^*(n)$ in the
basis of conjugacy classes of $\Gn$ are shown to be independent of
$n$, and in addition, they are non-negative integers \cite{Wa3}.
In a completely different approach \cite{LQW1}, the structure
constants of the cup product of the Heisenberg monomial basis for
$H^*(\Xn)$ are shown to be independent of $n$. In light of
Theorem~\ref{th:isom2}, these two ``independent of $n$" statements
are equivalent. Nevertheless, the positivity and integrality of
the structure constants in $R(\Gn)$, when transferred over to
$\Hn$, are still waiting for a geometric interpretation (even when
$\G$ is trivial).
\end{remark}

\section{\bf Moduli spaces of sheaves and fermionic Fock space}
\label{sec:moduli}
\subsection{The moduli spaces $\mathcal M (p,n)$ of sheaves}

Denote by $P$ the rank-$r$ integral lattice in $\mathbb H_1
=H^2_T(X)$ which is $\Z$-spanned by $\diamondsuit_i, 0\le i \le
r-1$, with the inherited bilinear form such that $\langle
\diamondsuit_i, \diamondsuit_j\rangle =\delta_{ij}$. We denote by
$\mathcal O_m$, $m \in\Z$, the line bundle $X \times_T \C$ over
$X$, where $T$ acts on $\C$ according to the representation
$\theta^m$.

\begin{lemma}
Given $p \in P$, there exists a $T$-equivariant line bundle over
$X$ whose $T$-equivariant first Chern class is $p$.
\end{lemma}

\begin{proof}
It suffices to show that there exists a $T$-equivariant line
bundle over $X$ whose $T$-equivariant first Chern class is
$\diamondsuit_i, 0 \le i \le r-1$. By (\ref{eq:chern_Lkn}),
\begin{eqnarray} \label{eq:first}
c_1^T(L_k) =\frac1r \left(\sum_{i=0}^{k-1} (k-r)
\diamondsuit_i +\sum_{i= k}^{r-1} k \diamondsuit_i \right ),
 \qquad 0 \le k \le r-1.
\end{eqnarray}

Combined with Lemma~\ref{H1}, this implies that
\begin{eqnarray} \label{eq:firstchern}
\diamondsuit_k
 &=& c_1^T(L_k)- c_1^T(L_{k+1}) +t,  \qquad 0 \le k< r-1 \\
\diamondsuit_{r-1}
 &=& c_1^T(L_{r-1}) +t. \nonumber
 \end{eqnarray}
Note that the $T$-equivariant first Chern class of the line bundle
$\mathcal O_m$ equals $mt$. Therefore, for $0 \le k< r-1$,
the $T$-equivariant first Chern class of
$L_k \otimes L_{k+1}^\vee \otimes \mathcal O_1$ is $\diamondsuit_k$,
and the $T$-equivariant first Chern class of $L_{r-1} \otimes
\mathcal O_1$ is $\diamondsuit_{r-1}$.
\end{proof}

\begin{remark}
We derive from Lemma~\ref{H1} and (\ref{eq:firstchern}) that
$$
([\Sigma_1], [\Sigma_2], \cdots, [\Sigma_{r-1}]) = - (c_1^T(L_1),
c_1^T(L_2), \cdots, c_1^T(L_{r-1})) \cdot C_{r-1}
$$
where $C_{r-1}$ is the Cartan matrix of type $A_{r-1}$. This can
be regarded as a equivariant geometric version of the McKay
correspondence \cite{McK}. A relation of this sort in the setup of
ordinary cohomology theory of the minimal resolutions was due to
Gonzalez-Sprinberg and Verdier \cite{GSV}.
\end{remark}

Denote by $\mathcal O(p)$ the $T$-equivariant line bundle
(guaranteed by the above Lemma) over $X$
whose $T$-equivariant first Chern class is $p \in P$.
We introduce the moduli space, denoted by $\mathcal M(p,n)$,
which parameterizes all rank-1 subsheaves of
$\mathcal O(p)$ such that the quotients are supported at finitely
many points of $X$ and have length $n$. Given $\mathcal I \in
\Xn$, then $\mathcal O(p) \otimes \mathcal I$ is an element in
$\mathcal M(p,n).$ Therefore, these moduli spaces for a fixed $n$
and varied $p$ are all $T$-equivariantly isomorphic to $\Xn$. In
particular, $\mathcal M(0,n) =\Xn$.
As before, the consideration of the equivariant cohomology ring
$H^*_T(\mathcal M(p, n))$ leads to a ring $\Hn^{(p)}
\stackrel{\text{def}}{=} H^{2n}_T(\mathcal M(p, n))$ whose product
is again denoted by $\star$. For $p \in P$, we introduce
\begin{eqnarray*}
\mathbb H^{(p)} =  \bigoplus_{n=0}^{\infty} \Hn^{(p)}, \qquad
 \FFock = \bigoplus_{p \in P} \mathbb H^{(p)}.
\end{eqnarray*}
The natural identification $\mathcal M(p, n) \cong \Xn$ leads to
the natural identification of the rings $\Hn^{(p)} \cong \Hn,$
which induces a bilinear form $\langle -, -\rangle^{(p)}$ on
$\Hn^{(p)}$ and $\mathbb H^{(p)}$. Similarly, we shall add the
superscript ${(p)}$ to denote the counterparts in $\mathcal M(p,
n)$ of the objects in $\Xn$ and its equivariant cohomology, such
as $\xi_\la^{(p)}$, $[\la]^{(p)}$, etc.

Given $\alpha \in P$, we denote by $S_\alpha$ the isomorphism from
$\Hn^{(p)}$ to $\Hn^{(p+\alpha)}.$ This induces isomorphisms
$S_\alpha: \FFock \rightarrow \FFock$ and $S_\alpha: \mathbb
H^{(p)} \rightarrow \mathbb H^{(p+\alpha)}$ for all $p\in P$.
Apparently, we have $S_{\alpha+\beta} =S_\alpha S_\beta$ for
$\alpha, \beta \in P$.
\subsection{The Chern character operators on fermionic Fock space}

Let $p \in P$.
Denote by $\pi_1$ and $\pi_2$ the projections of $\mathcal M(p,n)
\times X$ to the two factors. Over $\mathcal M(p,n)
\times X$, we have a universal exact sequence:
$$0 \longrightarrow \mathcal J(p) \longrightarrow \pi_2^*\;\mathcal
O(p) \longrightarrow  \mathcal Q(p) \to 0.$$
For $0\le k \le r-1$, we denote by $L_k(p)^{[n]}$ the
$T$-equivariant rank $n$ vector bundle over $\mathcal M(p,n)$
given by the pushforward $\pi_{1*}(\mathcal Q(p) \otimes \pi_2^*
L_k)$, whose fiber over a fixed point $\xi_\la^{(p)} \in \mathcal
M(p,n) \cong X^{[n]}$ is given by
\begin{eqnarray} \label{eq:fiber}
L_k(p)^{[n]} \mid_{\xi_\la^{(p)}} =\mathcal O(p) \otimes L_k^{[n]}
\mid_{\xi_\la^{(p)}}, \qquad \la \in \Pnr.
\end{eqnarray}

In the same way as defining the operators $\mathfrak G_k(z)$ and
$\mathfrak G_{k;m}$ which act on $\Hx =\mathbb H^{(0)}$, where $0
\le k \le r-1$ and $m \ge 0$, we can define the operators
$\mathfrak G_k^{(p)}(z)$ (respectively, $\mathfrak G_{k;m}^{(p)}$)
acting on $\mathbb H^{(p)}$ by the $\star$-product with $\sum_{m
\ge 0} t^{n-m} \ch^T_m (L_k(p)^{[n]}) \, z^m$ (respectively, with
$t^{n-m} \ch^T_m (L_k(p)^{[n]})$) on $\Hn^{(p)}$ for each $n$.

Take $p =\sum_i n_i \diamondsuit_i \in P$, where $n_i \in \Z$. One
can show that  as $T$-modules
$$\mathcal O(\diamondsuit_i)
|_{\xi_j} \cong \theta^{r\delta_{i,j} }.$$
Similarly as in Subsection~\ref{subsec:chern}, by the projection
formula and (\ref{eq:fiber}), we obtain
\begin{eqnarray} \label{eq:proj}
& & \ch_m^T(L_k(p)^{[n]}) \cup [\xi_\la^{(p)}] =  \frac{1}{m!}
 \left ( \sum_{i=0}^{k-1} \,  \sum_{\square \in \la^i}
   \big ( (n_ir +k-r+rc_{\square})t \big )^m \right.  \\
&& \qquad \qquad \qquad \qquad \qquad \left.
   + \sum_{i=k}^{r-1} \,\, \sum_{\square \in \la^i} \big (
   (n_ir +k+rc_{\square})t \big )^m  \right ) [\xi_\la^{(p)}].\nonumber
\end{eqnarray}

\begin{lemma}  \label{lem:sector}
Let $p =\sum_i n_i \diamondsuit_i \in P$, where $n_i \in \Z$, and
let $\la =(\la^0, \cdots, \la^{r-1}) \in \Pnr$ with $\la^i
=(\la_1^i, \la_2^i, \ldots)$. Then, for $0 \le k \le r-1$, we have
\begin{eqnarray*}
\mathfrak G_k^{(p)} (z) \cdot [\la]^{(p)}
 &=& \frac1{\varsigma (rz)}
  \left (
  \sum_{i=0}^{k-1}
  e^{(n_ir+k-r)z} \left ( \sum_{j=1}^\infty
     e^{(\la^i_j -j+1/2) rz} - \frac1{\varsigma (rz)} \right )
     \right.   \nonumber \\
 & &   \left.
   \qquad \quad + \sum_{i=k}^{r-1} e^{(n_ir+k)z} \left
(\sum_{j=1}^\infty
     e^{(\la^i_j -j+1/2) rz} - \frac1{\varsigma (rz)} \right )
     \right )
\; [\la]^{(p)}.
\end{eqnarray*}
\end{lemma}

\begin{proof}
Follows from (\ref{eq:alg}) and (\ref{eq:proj}) directly.
\end{proof}

For later purpose, we introduce the following modified classes
\begin{eqnarray}  \label{eq:modify}
{\ch}_{k;m}^{[n] (p)}
= t^{n-m}\ch^T_m(L_k(p)^{[n]}) +c^{(p)}_{k;m} t^n, \qquad m \ge -1
\end{eqnarray}
where the constant $c^{(p)}_{k;m}$ is defined by
$$ \sum_{m \ge -1} c^{(p)}_{k;m} z^m
 = \frac1{\varsigma(rz)^2}
  \left( e^{(k-r)z}\sum_{i=0}^{k-1} (e^{n_irz}-1)
  + e^{kz} \sum_{i=k}^{r-1} (e^{n_irz}-1) \right )$$
and $\ch^T_{-1}(L_k(p)^{[n]})=0$ by convention. Equivalently, if
we define
\begin{align}  \label{eq:tilde}
\widetilde{\mathfrak G}^{(p)}_k(z) =& \mathfrak G^{(p)}_k(z) + \\
& \frac1{\varsigma(rz)^2}
  \left( e^{(k-r)z}\sum_{i=0}^{k-1} (e^{n_irz}-1)
  + e^{kz} \sum_{i=k}^{r-1} (e^{n_irz}-1) \right ) \text{Id}
  \nonumber
\end{align}
and further write
$$\widetilde{\mathfrak G}^{(p)}_k(z) =
\sum_{m = -1}^\infty \widetilde{\mathfrak G}^{(p)}_{k;m} z^m, $$
then
$$\widetilde{\mathfrak G}^{(p)}_{k;m} ={\mathfrak G}^{(p)}_{k;m}
+c^{(p)}_{k;m} \text{Id}, $$ (where we take the convention that
${\mathfrak G}^{(p)}_{k;-1}=0$),
and $\widetilde{\mathfrak G}^{(p)}_{k;m}$ acts on $\Hn^{(p)}$ by
the product with ${\ch}_{k;m}^{[n] (p)}$. When $p=0$, we have
$$c^{(0)}_{k;m} =0, \quad {\ch}_{k;m}^{[n] (0)}={\ch}_{k;m}^{[n]},
\quad \widetilde{\mathfrak G}^{(0)}_{k;m}  = {\mathfrak G}_{k;m},
\quad \widetilde{\mathfrak G}^{(0)}_k(z) = {\mathfrak G}_k(z). $$

By the standard boson-fermion correspondence (cf. \cite{MJD}), we
shall identify $\FFock$ with the fermionic Fock space of $r$ pairs
of fermions $\psi^{+,k}(z),\psi^{-,k}(z), 0\le k \le r-1$. More
explicitly, we extend the Heisenberg operators $\mathfrak
p_n(\alpha)$, where $\alpha \in \mathbb H_1$ and $n \neq 0$, to
act on $\mathbb H^{(p)}, p\in P$ via the identification $\mathbb
H^{(p)} \cong \Hx$, and let $\mathfrak p_0 (\alpha)$ act on
$\mathbb H^{(p)}$ by $\langle \alpha, p \rangle \, \text{Id}$.
Note that this is compatible with $\mathfrak p_0 (\alpha)=0$ on
$\Hx$ thanks to $ \Hx=\mathbb H^{(0)}$. The fermionic fields can
be expressed as
$$ \psi^{\pm ,k}(z)
 = S_{\diamondsuit_k} z^{\pm \mathfrak p_0(\diamondsuit_k)}
 \exp \left ( \sum_{n < 0} \mp \frac{\mathfrak p_n(\diamondsuit_k)}{n} z^{-n}
  \right )
  \exp \left ( \sum_{n > 0} \mp \frac{\mathfrak p_n(\diamondsuit_k)}{n} z^{-n}
  \right ).
  $$

It has been well known (cf. e.g. \cite{MJD}) that the completed
infinite-rank general linear Lie algebra $\widehat{gl}_\infty$,
whose standard basis is denoted by $E_{i,j}$, $i,j\in\Z +1/2$,
acts on the fermionic Fock space of a pair of fermions. Thus, we
have an {\em orthogonal} direct sum of $r$ copies of
$\widehat{gl}_\infty$, denoted by $\widehat{gl}_\infty^{r}$,
acting on $\FFock$, whose generators will be denoted by
$E_{i,j}^{(k)}, 0\le k \le r-1$. Using the generating field
$$ E^{(k)} (z,w)
 = \sum_{i,j\in\Z +1/2} E_{i,j}^{(k)} z^{i -1/2} w^{-j-1/2}$$
the action of $\widehat{gl}^{r}_\infty$ is simply given by
 $$ E^{(k)} (z,w) = \; :\psi^{+,k}(z) \psi^{-,k}(w):.$$
We refer to {\em loc. cit.} for the standard notation of normal
ordering $:\;\;:$. Denote
$$ \mathcal E^{(i)}(z) = \frac1{\varsigma (rz)} \sum_{m \in \Z
+\frac12} e^{mrz} E_{m,m}^{(i)}.$$
We introduce the following operator in $\End (\FFock)$:
\begin{eqnarray}  \label{H_k}
\mathfrak H_{k} (z)  =
   e^{(k-r)z}\sum_{i=0}^{k-1}\mathcal E^{(i)}(z)
      + e^{kz} \sum_{i=k}^{r-1} \mathcal E^{(i)}(z)
\end{eqnarray}
which is further expanded as
\begin{eqnarray*}
\mathfrak H_{k} (z) =
\sum_{m=-1}^\infty \mathfrak H_{k;m} z^m,\qquad 0 \le k \le r-1.
\end{eqnarray*}
Alternatively, we can rewrite (\ref{H_k}) as
\begin{eqnarray}  \label{H_k2}
\mathfrak H_{k} (z)  =
   \sum_{i=0}^{r-1} e^{a(k,i)z} \mathcal E^{(i)}(z)
\end{eqnarray}
where
\begin{eqnarray} \label{eq:number}
a(k,i) = \left \{
 \begin{array}{ll}
  k-r, & \quad \text{ if } 0 \le i \le k-1  \\
  k, &  \quad \text{ if }  k\le i \le r-1
 \end{array}
 \right.
\end{eqnarray}
 For example, $\mathfrak H_{k;-1}
=\sum_{i=0}^{r-1} r^{-1} \mathfrak p_0(\diamondsuit_i).$ It is
well known (cf. \cite{MJD}) that the elements in the Cartan
subalgebra of $\widehat{gl}_\infty$ diagonalize the basis of Schur
functions. Thanks to the identification of the $[\la]$'s with
multi-variable Schur functions in Theorem~\ref{th:isom1}, we have
the following straightforward multi-variable generalization in our
context (cf. Lemmas 3.1 and 3.6 in \cite{LQW2}; also cf.
\cite{OP}).

\begin{lemma}   \label{lem:eigen}
The operators $\,\mathcal E^{(i)}(z)$ and $\mathfrak H_{k} (z)$
diagonalize the elements $[\la]$, $\la \in \cup_n \mathcal
P_n(r)$. More explicitly, for $\la =(\la^0, \cdots, \la^{r-1})$
with $\la^i =(\la^i_1, \la^i_2, \cdots)$, we have
 $$ \mathcal E^{(i)}(z) \cdot [\la]
 =  \frac1{\varsigma (rz)}\left (\sum_{j=1}^\infty
     e^{(\la^i_j -j+1/2) rz} - \frac1{\varsigma (rz)}
     \right) \;[\la]. \qquad\qquad\qquad \qed
 $$
\end{lemma}

\begin{theorem} \label{th:G=H}
 Given $p\in P$, we have the identification
$$\mathfrak H_k(z) \mid_{\mathbb H^{(p)}} =
\widetilde{\mathfrak G}_k^{(p)} (z).$$
\end{theorem}

\begin{proof}

First consider $p=0$. It follows from
Proposition~\ref{prop:chern}, the definition (\ref{H_k}) of
$\mathfrak H_k (z)$, and Lemma~\ref{lem:eigen} that $\mathfrak
H^{(k)}(z)$ and $\mathfrak G_k(z)$ have the same eigenvalues on
the basis elements $[\la]$ of $\mathbb H^{(0)}$, whence
$$\mathfrak H^{(k)}(z) \mid_{\mathbb H^{(0)}} =\mathfrak G_k (z) =
\widetilde{\mathfrak G}_k^{(0)} (z).$$

Now let $p =\sum_i n_i \diamondsuit_i$, $n_i \in\Z$. For $d
\in\Z$, we can show by a standard computation that (cf. Lemma~3.5,
\cite{LQW2})
\begin{eqnarray*}
S_{\diamondsuit_i}^{-d} \, \mathcal E^{(j)}(z) \,
S_{\diamondsuit_i}^d
&=& \left \{
 \begin{array}{ll}
 e^{drz}\mathcal E^{(i)}(z) +
 \frac{e^{drz}-1}{\varsigma(rz)^2} \; \text{\rm Id}, & \quad \text{ if
}\; i=j  \\
  \mathcal E^{(j)}(z), & \quad \text{ if } \; i \neq j
 \end{array}
 \right.
 \end{eqnarray*}
Thus, using (\ref{H_k}) we have
\begin{eqnarray*}
S_{p}^{-1} \, \mathfrak H_k (z) \, S_{p} %
&=& \prod_i S_{\diamondsuit_i}^{-n_i} \, \mathfrak H_k (z) \prod_i
S_{\diamondsuit_i}^{n_i}  \\
&=&
   e^{(k-r)z} \sum_{i=0}^{k-1} e^{n_irz} \mathcal E^{(i)}(z)
      + e^{kz}  \sum_{i=k}^{r-1} e^{n_irz} \mathcal E^{(i)}(z) \\
& & +  \left (\sum_{i=0}^{k-1} \frac{e^{(k-r)z} (e^{n_i r
z}-1)}{\varsigma(rz)^2}+ \sum_{i=k}^{r-1}  \frac{e^{kz}(e^{n_ir
z}-1)}{\varsigma(rz)^2} \right) \text{\rm Id}.
\end{eqnarray*}
Together with Lemma~\ref{lem:eigen}, this implies by a little
algebra manipulation that
\begin{eqnarray*}
\mathfrak H_k (z)\cdot [\la]^{(p)}
&=& \mathfrak H_k (z) S_{p} \cdot [\la] \\
&=& \frac1{\varsigma (rz)}
  \left (
  \sum_{i=0}^{k-1}
   \left ( \sum_{j=1}^\infty e^{(n_ir+k-r)z}
     e^{(\la^i_j -j+1/2) rz} - \frac{e^{(k-r)z}}{{\varsigma (rz)}}
\right )
     \right.   \nonumber \\
 & &   \left.
   \qquad \quad + \sum_{i=k}^{r-1} \left (\sum_{j=1}^\infty
     e^{(n_ir+k)z}e^{(\la^i_j -j+1/2) rz}
     - \frac{e^{kz}}{\varsigma (rz)} \right )
     \right )
\; [\la]^{(p)}.\end{eqnarray*}

On the other hand, by (\ref{eq:tilde}) and Lemma~\ref{lem:sector},
we have
\begin{eqnarray*}
\widetilde{\mathfrak G}_k^{(p)} (z) \cdot [\la]^{(p)}
 &=& \frac1{\varsigma (rz)}
  \left (
  \sum_{i=0}^{k-1}
   \left ( \sum_{j=1}^\infty e^{(n_ir+k-r)z}
     e^{(\la^i_j -j+1/2) rz} - \frac{e^{(k-r)z}}{\varsigma (rz)} \right
)
     \right.   \nonumber \\
 & &   \left.
   \qquad \quad + \sum_{i=k}^{r-1} \left (\sum_{j=1}^\infty
     e^{(n_ir+k)z}e^{(\la^i_j -j+1/2) rz}
     - \frac{e^{kz}}{\varsigma (rz)} \right )
     \right )
\; [\la]^{(p)}.
\end{eqnarray*}
Thus, we see that the two operators $\mathfrak H_k (z)$ and
$\widetilde{\mathfrak G}_k^{(p)} (z)$ have the same eigenvalues on
every $[\la]^{(p)}$. This proves the theorem.
\end{proof}
\section{\bf Generating functions of intersection numbers}
\label{sec:inter}
\subsection{The $N$-point functions of intersection numbers}

We are interested in the equivariant intersection numbers for
$\Xn$ (associated to $\la, \mu \in \Pnr$):
\begin{eqnarray*}
\left \langle \la, {\ch}_{k_1;m_1}^{[n]} \cdots
 {\ch}_{k_N;m_N}^{[n]},\mu \right \rangle
 : =
\left\langle \mathfrak p_{-\la}, {\ch}_{k_1;m_1}^{[n]} \star
\cdots \star {\ch}_{k_N;m_N}^{[n]}\star \mathfrak p_{-\mu}
\right\rangle
\end{eqnarray*}
where the classes ${\ch}_{k_i;m_i}^{[n]}$ are defined in
(\ref{ch_km_n}). We can organize these intersection numbers
into the {\em $N$-point function}
(for fixed $0 \le k_1, \ldots, k_N \le r-1$):
\begin{eqnarray*}
& &G_{\la, \mu;r} (z_1,\ldots, z_N;k_1,\ldots, k_N)  \\
&=&\sum_{m_1,\ldots, m_N =0}^\infty z_1^{m_1}\cdots z_N^{m_N}
   \left \langle \la,{\ch}_{k_1;m_1}^{[n]} \cdots
   {\ch}_{k_N;m_N}^{[n]},\mu \right \rangle,
\end{eqnarray*}
which can be further reformulated in an operator
form thanks to Theorem~\ref{th:G=H}:
\begin{eqnarray}  \label{eq:opform}
G_{\la, \mu;r} (z_1,\ldots, z_N;k_1,\ldots, k_N)
 &=& \left \langle \mathfrak p_{-\la},\mathfrak
G_{k_1}(z_1) \cdots \mathfrak G_{k_N} (z_N)
\mathfrak p_{-\mu} \right \rangle  \\
&=& \langle \mathfrak p_{-\la},\mathfrak H_{k_1}(z_1) \cdots
\mathfrak H_{k_N} (z_N) \mathfrak p_{-\mu} \rangle.  \nonumber
\end{eqnarray}
For latter purpose, we have indicated $r$ explicitly in the
notation above. When $r=1$ (i.e. $X=\C^2$), we simply write it as
$G_{\la, \mu;1} (z_1,\ldots, z_N)$.

We shall need the notion of a partition of the set $\{1, \ldots,
N\}$ which is a collection of subsets $\mathcal N :=(\mathcal N_0,
\ldots, \mathcal N_{r-1})$ whose disjoint union is $\{1, \ldots,
N\}.$ Denote ${\bf z} =(z_1, \cdots, z_N)$ and ${\bf k} =(k_1,
\cdots, k_N)$, and let
$$F({\bf z}, {\bf k}; \mathcal N) :=\prod_{i=0}^{r-1}
\prod_{j\in\mathcal N_i} e^{a(k_j,i)z_j}.$$
(See (\ref{eq:number}) for notations). The next proposition
reduces the calculation of the $N$-point function for a general
$r$ to the case $r=1$.

\begin{proposition} \label{prop:reduct}
Fix $\la=(\la^0,\ldots, \la^{r-1})$ and  $\mu=(\mu^0,\ldots,
\mu^{r-1})$ in $\Pnr$. Then
\begin{eqnarray*}
G_{\la, \mu;r} (z_1,\ldots, z_N;k_1,\ldots, k_N)=0
\end{eqnarray*}
unless $|\la^i| =|\mu^i|$ for $0 \le i \le r-1$; in this case, then
\begin{eqnarray*}
G_{\la, \mu;r} (z_1,\ldots, z_N;k_1,\ldots, k_N)
&=& \sum_{\mathcal N} \left ( F({\bf z}, {\bf k}; \mathcal N)
\prod_{i=0}^{r-1}
   G_{\la^i, \mu^i;1} (r z_{l_{i,1}},\ldots, r z_{l_{i,n_i}})
   \right)
\end{eqnarray*}
summed over all partitions $\mathcal N=(\mathcal N_0, \ldots,
\mathcal N_{r-1})$ of the set $\{1, \ldots, N\}$, where $\mathcal
N_i =\{l_{i,1}, \ldots, l_{i,n_i} \}$.
\end{proposition}

\begin{proof}
We can write $\mathfrak a^R_{-\mu}$ as a linear combination:
$\mathfrak a^R_{-\mu} =\sum_{\nu \in \Pnr} d_{\mu,\nu}s_\nu,$
where $\nu$ satisfies $|\nu^i|=|\mu^i|$ for $0 \le i \le r-1$, cf.
Section~\ref{sec:wreath}. By Theorem~\ref{th:isom1}, we have
$\mathfrak p_{-\mu} =\sum_{\nu \in \Pnr} d_{\mu,\nu} [\nu]$. The
same can be done for $\mathfrak p_{-\la}$. Since the $[\nu]$'s are
orthogonal eigenvectors for the operator $\mathfrak G_{k} (z)$, we
see that
$$G_{\la, \mu;r} (z_1,\ldots, z_N;k_1,\ldots, k_N)=0$$
unless $|\la^i| =|\mu^i|$ for $0 \le i \le r-1$.

Now assume $|\la^i| =|\mu^i|$ for $0 \le i \le r-1$. Recall from
(\ref{p-la}) that $$\mathfrak p_{-\mu} = \prod_i \mathfrak
p_{-\mu^i}(\diamondsuit_i) \vac.$$ Note that the operator
$\mathcal E^{(j)}(z)$ commutes with the operators
$\mathcal E^{(i)}(z)$ and $ \mathfrak
p_{-\mu^i}(\diamondsuit_i)$ for $j \neq i$. Thus, by (\ref{H_k2})
and (\ref{eq:opform}), we obtain
\begin{align} \label{eq:step1}
 G_{\la, \mu;r} &(z_1,\ldots, z_N;k_1,\ldots, k_N)   \\
 &= \sum_{\mathcal N}
 \left \langle \mathfrak p_{-\la},
 \prod_{i=0}^{r-1} \prod_{j\in\mathcal N_i}
 \left( e^{a(k_j,i)z_j} \mathcal E^{(i)}(z_j)\right)
 \mathfrak p_{-\mu}
 \right \rangle  \nonumber \\
 &= \sum_{\mathcal N}
 F({\bf z}, {\bf k}; \mathcal N) \cdot
 \left \langle   \mathfrak p_{-\la},
 \prod_{i=0}^{r-1}  \prod_{j\in\mathcal N_i}
 \mathcal E^{(i)}(z_j)
 \mathfrak p_{-\mu^i}(\diamondsuit_i)\vac
 \right \rangle \nonumber
\end{align}
summed over all partitions $\mathcal N =(\mathcal N_0, \ldots,
\mathcal N_{r-1})$ of the set $\{1, \ldots, N\}.$

Observe that we can write $\prod_{j\in\mathcal N_i} \mathcal
E^{(i)}(z_j) \mathfrak p_{-\mu^i}(\diamondsuit_i)\vac$ as a linear
combination of $\mathfrak p_{-\nu^i}(\diamondsuit_i) \vac$'s, and
thus, $\prod_{i=0}^{r-1} \prod_{j\in\mathcal N_i} \mathcal
E^{(i)}(z_j) \mathfrak p_{-\mu^i}(\diamondsuit_i)\vac$ as a linear
combination of $\mathfrak p_{-\nu}$'s with $\nu=(\nu^0,\ldots,
\nu^{r-1}) \in \Pnr$. Also observe that
$$ \left\langle  \mathfrak p_{-\la},
 \mathfrak p_{-\nu} \right\rangle
 =\prod_i \big \langle \mathfrak p_{-\la^i}(\diamondsuit_i) \vac,
\mathfrak p_{-\nu^i}(\diamondsuit_i) \vac \big \rangle.
$$
It follows from these observations and (\ref{eq:step1}) that
\begin{align*}
G_{\la, \mu;r}& (z_1,\ldots, z_N;k_1,\ldots, k_N) \\
 &= \sum_{\mathcal N} F({\bf z}, {\bf k}; \mathcal N) \cdot
 \prod_{i=0}^{r-1}
  \left \langle \mathfrak p_{-\la^i}(\diamondsuit_i)\vac,
  \prod_{j\in\mathcal N_i} \mathcal E^{(i)}(z_j)
   \mathfrak p_{-\mu^i}(\diamondsuit_i)
 \vac \right \rangle  \\
&= \sum_{\mathcal N} \left(F({\bf z}, {\bf k}; \mathcal N) \cdot
\prod_{i=0}^{r-1} G_{\la^i, \mu^i;1} (r z_{l_{i,1}},\ldots, r
z_{l_{i,n_i}}) \right).
\end{align*}
This completes the proof.
\end{proof}

For $N=1$ the above proposition reads
\begin{eqnarray} \label{eq:N=1}
G_{\la, \mu;r} (z;k)
&=& \sum_{i=0}^{r-1} e^{z a(k,i)}
   G_{\la^i, \mu^i;1} (r z).
\end{eqnarray}

\begin{remark}
When $r=1$, the $N$-point function $G_{\la^i, \mu^i;1} (r
z_1,\ldots, r z_N)$ of intersection numbers on Hilbert schemes of
points on the affine plane has been studied in some detail in
\cite{LQW2}. It was shown in {\em loc. cit.} that the $N$-point
functions for the affine plane has a precise connection with the
$N$-pointed functions of disconnected Gromov-Witten invariants of
$\mathbb P^1$ studied in \cite{OP}. In addition, an explicit
formula for the $1$-point function was given in Theorem~4.2,
\cite{LQW2}. When combined with (\ref{eq:N=1}), we obtain the
$1$-point function for a general $r$.
\end{remark}

\subsection{The $\tau$-functions of intersection numbers on $\mathcal
M(p,n)$}

For $0 \le i \le r-1$, we introduce two sequences of indeterminates:
\begin{eqnarray*}
t_i= (t_{i,1}, t_{i,2}, \ldots), \qquad
s_i =(s_{i,1}, s_{i,2}, \ldots).
\end{eqnarray*}
Set $t =(t_0,\ldots, t_{r-1}), s=(s_0,\ldots,
s_{r-1})$. Define the following half vertex operators:
$$\Gamma^{(i)}_\pm (t) = \exp \left( \sum_{m>0} \frac1m t_{i,m}
\mathfrak p_{\pm m}(\diamondsuit_i) \right).$$

Define the numbers  $n_{k;d}$, where $d \ge -1, 0 \le k \le r-1$,
by the generating function
\begin{eqnarray}  \label{add1}
\sum_{d=-1}^\infty n_{k;d} \, z^d = \frac{e^{-kz}}{1-e^{-rz}}.
\end{eqnarray}
Introduce the following elements in $\mathbb H_n^{(p)}$, which are
linear combinations of equivariant Chern characters of
$L_k(p)^{[n]}$ (for $0\le k \le r-1$ and $m \ge -1$):
\begin{eqnarray}  \label{add2}
 %
%
\widetilde{\ch}^{[n](p)}_{k;m}
&= & \sum_{d=-1}^{m+1}  \left( n_{k;d} \, \ch^{[n](p)}_{k;m-d} -
 n_{k+1;d} \, \ch^{[n](p)}_{k+1;m-d}\right )
\end{eqnarray}
where we adopt the convention $\ch^{[n](p)}_{r;m-d}
= \ch^{[n](p)}_{0;m-d}$.
As we shall see in Proposition~\ref{prop:geom}, these modified
Chern characters admit explicit geometric constructions.

Given a multi-partition $\mu =(\mu^0, \cdots, \mu^{r-1})$ with
$\mu^i =(\mu^i_1, \mu^i_2, \ldots)$, we define
\begin{eqnarray*}
t^{(i)}_{\mu_i} = t_{i,\mu^i_1} t_{i,\mu^i_2}\cdots, \qquad
t_\mu =\prod_i t^{(i)}_{\mu_i},
\end{eqnarray*}
and similarly define $s_\mu$. Let $x_i=(x_{i,0},
x_{i,1}, x_{i,2}, \ldots)$, $0 \le i \le r-1$ and $x=(x_0,\ldots,
x_{r-1})$, be some other sequences of indeterminates. We introduce
the following generating function for the equivariant intersection
numbers on $\mathcal M(p,n)$:
\begin{eqnarray*}
\tau (x,t,s,p)
 =\sum_n \sum_{ \|\la\| =\|\mu\| =n}
 t_\la s_\mu \left\langle \la, \exp
 \left(\sum_{k=0}^{r-1}\sum_{m=0}^\infty x_{k,m}
\widetilde{\ch}_{k;m}^{[n](p)} \right), \mu \right\rangle^{(p)}
\end{eqnarray*}

By solving (\ref{H_k}) we obtain for $0\le k \le r-1$ that
\begin{eqnarray*}
%
\mathcal E^{(k)}(z)
 &=& \frac{1}{1 -e^{-rz}} \left(e^{-kz} \mathfrak H_{k}(z)
 -e^{-(k+1)z}\mathfrak H_{k+1}(z) \right)
\end{eqnarray*}
where we let $\mathfrak H_r(z) = \mathfrak H_0(z)$.
Denote by $\mathcal E^{(k)}(z) =\sum_{m=-1}^\infty \mathcal
E_{k;m}z^m$. Then,
\begin{eqnarray*}
   \mathcal E_{k;m} (t^{n(p)})
&=&\text{Coeff}_{z^m} \left \{ \mathcal E^{(k)}(z)(t^{n(p)})
   \right \} \\
&=&\text{Coeff}_{z^m} \left \{ \frac{1}{1 -e^{-rz}}
   \left(e^{-kz} \mathfrak H_{k}(z)(t^{n(p)})
   -e^{-(k+1)z}\mathfrak H_{k+1}(z)(t^{n(p)}) \right) \right \}
\end{eqnarray*}
where $t^{n(p)}$ denotes
the image of $t^n$ under the isomorphism $S_p: \mathbb H_n
\rightarrow \mathbb H_n^{(p)}$. Combining this with
Theorem~\ref{th:G=H}, we obtain
\begin{eqnarray*}
  \mathcal E_{k;m} (t^{n(p)})
= \text{Coeff}_{z^m} \left \{ \frac{e^{-kz}}{1 -e^{-rz}}
   \widetilde{\mathfrak G}_k^{(p)}(z)(t^{n(p)})
   -\frac{e^{-(k+1)z}}{1 -e^{-rz}}
\widetilde{\mathfrak G}_{k+1}^{(p)} (z)(t^{n(p)}) \right \}.
\end{eqnarray*}
By the definition of the operator $\widetilde{\mathfrak G}_k^{(p)}(z)$,
(\ref{add1}) and (\ref{add2}), we get
\begin{eqnarray}  \label{eq:iden}
\mathcal E_{k;m} (t^{n(p)}) = \widetilde{\ch}^{[n](p)}_{k;m}.
\end{eqnarray}

\begin{proposition}  \label{prop:geom}
For $\alpha \in H^*_T(X)$, let $G^T_m (\alpha, n)$ be the
$H^{2m}_T(\Xn)$-component of the class
$\pi_{1*} \left(\ch^T(\mathcal O_{\mathcal Z_n})\cup \pi_2^*\alpha
\cup \pi_2^*\td^T (X) \right)$, where $\td^T (X)$ is the
$T$-equivariant Todd class of $X$. Then,
$\widetilde{\ch}^{[n](0)}_{k;m}=t^{n-m-1} G^T_{m+1}
(\diamondsuit_k/r, n)$ for $m \ge -1$.
\end{proposition}
%
%
\begin{proof}
First, applying the equivariant Grothendieck-Riemann-Roch Theorem
\cite{EGr} to $L_k^{[n]} = \pi_{1*}(\mathcal O_{\mathcal Z_n}
\otimes \pi_2^*L_k)$, we have
$$\ch^T (L_k^{[n]})
= \pi_{1*} \left(\ch^T(\mathcal O_{\mathcal Z_n})\cup
\pi_2^*\ch^T(L_k) \cup \pi_2^*\td^T (X) \right).
$$
It follows that
\begin{eqnarray} \label{GRR}
G_m^T (\ch^T(L_k), n) = \ch_m^T (L_k^{[n]}).
\end{eqnarray}

Next, recall that the operator $\mathfrak G_k(z) \in \End (\Hx)$
is defined by the $\star$-product with the class $\sum_{m\ge 0}
t^{n-m}\ch_m^T (L_k^{[n]}) z^m$ in $\mathbb H_n$ for every $n$.
We define the operator
\begin{eqnarray*}
\mathfrak G^{(\alpha)}(z) \in \End (\Hx)
\end{eqnarray*}
by the $\star$-product with $\sum_{\ell \ge 0} t^{n-\ell} G^T_\ell
(\alpha, n) z^{\ell-1}$ in $\mathbb H_n$ for every $n$.

Note that $\diamondsuit_i/r \cup \diamondsuit_j/r = \delta_{i, j}
\cdot t \,\, \diamondsuit_i/r$. By (\ref{eq:first}) and
Lemma~\ref{H1}~(i), we have
\begin{eqnarray} \label{char}
\ch^T(L_k)
= \sum_{i=0}^{k-1} \sum_{d=0}^\infty
\frac{(k-r)^d t^{d-1}}{d!} \diamondsuit_i/r
  +\sum_{i=k}^{r-1} \sum_{d=0}^\infty
\frac{k^d t^{d-1}}{d!} \diamondsuit_i/r.
\end{eqnarray}
Therefore, by (\ref{GRR}), (\ref{char}) and the linearity of
$G^T_m(\alpha,n)$ on $\alpha$, we have
\begin{eqnarray*}
 && \sum_{m=0}^\infty t^{n-m}\ch_m^T (L_k^{[n]}) z^m \\
 &=& \sum_{m=0}^\infty t^{n-m} G_m^T (\ch^T(L_k), n) z^m \\
 &=&  \sum_{m=0}^\infty z^m t^{n-m}\sum_{d=0}^\infty \left(\sum_{i=0}^{k-1}
  \frac{(k-r)^d }{d!} G_m^T (t^{d-1} \diamondsuit_i/r, n)
  + \sum_{i=k}^{r-1}
 \frac{k^d}{d!}G_m^T (t^{d-1} \diamondsuit_i/r, n)\right) \\
 &=& \sum_{m=0}^\infty z^m t^{n-m}\sum_{d=0}^m \left(\sum_{i=0}^{k-1}
  \frac{(k-r)^d t^{d-1} }{d!} G_{m-d+1}^T (\diamondsuit_i/r,
 n)
  \right. \\
 &&
 \left. \qquad\qquad\qquad\qquad\qquad
  + \sum_{i=k}^{r-1}
 \frac{k^dt^{d-1}}{d!}G_{m-d+1}^T (\diamondsuit_i/r, n)\right)
\end{eqnarray*}
which, by setting $\ell = m-d+1$, is equal to
\begin{eqnarray*}
& &\sum_{i=0}^{k-1} \sum_{d=0}^\infty
   \frac{(k-r)^dz^d}{d!} \sum_{\ell =0}^\infty z^{\ell-1}
   t^{n -\ell} G_\ell^T (\diamondsuit_i/r, n)    \\
&+&\sum_{i=k}^{r-1} \sum_{d=0}^\infty
   \frac{k^dz^d}{d!} \sum_{\ell =0}^\infty z^{\ell-1}
   t^{n -\ell} G_\ell^T (\diamondsuit_i/r, n).
\end{eqnarray*}
It follows from the definition of the operator $\mathfrak G$'s
that
\begin{eqnarray} \label{formula}
\mathfrak G_k(z) =\sum_{i=0}^{k-1} e^{(k-r)z} \mathfrak
G^{(\diamondsuit_i/r)}(z) +\sum_{i=k}^{r-1} e^{kz}\mathfrak
G^{(\diamondsuit_i/r)}(z).
\end{eqnarray}
Thus, $\mathfrak G^{(\diamondsuit_i/r)}(z)
=\mathcal E^{(i)}(z) \mid_{\Hx}$
by applying Theorem~\ref{th:G=H}, solving (\ref{H_k}) for
$\mathcal E^{(i)}(z)$ and solving (\ref{formula}) for
$\mathfrak G^{(\diamondsuit_i/r)}(z)$.
Combining this with (\ref{eq:iden}), we obtain
\begin{eqnarray*}
  \widetilde{\ch}^{[n](0)}_{k;m}
= \mathcal E_{k;m} (t^{n})
= \text{Coeff}_{z^m} \left \{ \mathcal E^{(k)}(z)(t^{n}) \right \}
= \text{Coeff}_{z^m} \left \{ \mathfrak G^{(\diamondsuit_i/r)}(z)
  (t^{n}) \right \}.
\end{eqnarray*}
It follows from the definition of $\mathfrak
G^{(\diamondsuit_i/r)}(z)$ that
$\widetilde{\ch}^{[n](0)}_{k;m}=t^{n-m-1} G^T_{m+1}
(\diamondsuit_k/r,n)$.
\end{proof}

\begin{remark} \label{rem:isom}
We can enhance the ring isomorphism $\phi: \Hn \longrightarrow
R(\Gn)$ in Theorem~\ref{th:isom1} as follows. By comparing
Proposition~\ref{prop:JM} and Lemma~\ref{lem:eigen}, we see that
$$\phi (\mathcal E^{(i)}(z) \cdot a) =
 \mathfrak O^{(\g_i)} (z) \cdot \phi (a).$$
It follows that
$$\phi \big(\widetilde{\ch}^{[n](0)}_{k;m} \big) =
\Xi^m_n(\g_k).$$
Thus, as a counterpart of a result in \cite{Wa2}, the classes
$\widetilde{\ch}^{[n](0)}_{k;m}$, where $0 \le m <n$ and $0 \le k
\le r-1$, form a set of ring generators for $\mathbb H_n$.
\end{remark}

\begin{theorem}
\begin{enumerate}
 \item[(i)]
Let $p \in P$. The $\tau$-function $\tau (x,t,s, p)$ affords the
following operator formulation:
\begin{eqnarray*}
  \left \langle \vac, S_p^{-1}
 \prod_{k=0}^{r-1} \Gamma^{(k)}_+(t) \cdot \exp \left(\sum_{k=0}^{r-1}
 \sum_{m=0}^\infty x_{k,m}
\mathcal E_{k;m} \right)\cdot \prod_{k=0}^{r-1}
\Gamma^{(k)}_-(s)\, S_p\, \vac \right\rangle
\end{eqnarray*}
\item[(ii)] Write $p =\sum_i n_i \diamondsuit_i$ with $n_i \in
\Z$. Then,
$$\tau (x,t,s,p) = \prod_{k=0}^{r-1} \tau (x_k,t_k,s_k, n_k)$$
Here $\tau (x_k,t_k,s_k, n_k)$ denotes
\begin{eqnarray*}
 \left \langle \vac,
S_{\diamondsuit_k}^{-{n_k}}
  \Gamma^{(k)}_+(t) \cdot \exp \left(\sum_{m=0}^\infty x_{k,m}
\mathcal E_{k;m} \right) \cdot \Gamma^{(k)}_-(s) \,
S_{\diamondsuit_k}^{n_k} \, \vac \right\rangle
\end{eqnarray*}
 and it satisfies the $2$-Toda hierarchy
of Ueno-Takasaki \cite{UT}.
\end{enumerate}
\end{theorem}

\begin{proof}
Note that
$$\Gamma^{(i)}_-(s) =\sum_{n \ge 0} \sum_{|\la^i|=n}
t^{(i)}_{\la^i} \mathfrak p_{-\la^i}(\diamondsuit_i)$$
and $\Gamma^{(i)}_+(t)$ is the adjoint operator of
$\Gamma^{(i)}_-(t)$. Now Part (i) follows from (\ref{eq:iden}).

Observe that the bilinear form has the factorization property:
$$\langle [\la], [\mu] \rangle
 = \prod_{k=0}^{r-1} \langle [\la^k], [\mu^k] \rangle$$ for
$\la =(\la^0, \cdots, \la^{r-1})$ and $\mu =(\mu^0, \cdots,
\mu^{r-1})$. Now Part (ii) follows from the fact that the
operators associated to different $k$ in part (i) commute.

Note that the functions $\tau (x_k,t^k,s^k, n_k)$ can be
interpreted as generating functions of the equivariant
intersection numbers on Hilbert schemes of points on the affine
plane (i.e. $r=1$ in our case) studied in \cite{LQW2}. It was
shown in {\em loc. cit.} that $\tau (x_k,t^k,s^k, n_k)$ satisfies
the $2$-Toda hierarchy of Ueno-Takasaki \cite{UT}.
\end{proof}


\begin{thebibliography}{ABCD}

\bibitem[EGr]{EGr} D. Edidin and W. Graham,
{\em Riemann-Roch for equivariant Chow groups}, Duke Math. J. {\bf
102} (2000), 567--594.

\bibitem[ES]{ES} G. Ellingsrud and S. Str\o mme,
{\em On the homology of the Hilbert scheme of points in the
plane}, Invent. Math. {\bf 87} (1987), 343--352.

\bibitem[EG]{EG} P. Etingof and V. Ginzburg,
{\em Symplectic reflection algebras, Calogero-Moser space, and
deformed Harish-Chandra homomorphism}, Invent. Math. {\bf 147}
(2002), 243--348.

\bibitem[FJW]{FJW} I.~Frenkel, N.~Jing and W.~Wang,
{\em Vertex representations via finite groups and the McKay
correspondence}, Internat. Math. Res. Notices, {\bf 4} (2000),
195--222.

\bibitem[FLM]{FLM} I. Frenkel, J. Lepowsky and A. Meurman,
{\em Vertex operator algebras and the Monster,} Academic Press,
New York, 1988.

\bibitem[GSV]{GSV} G.~Gonzalez-Sprinberg and J.-L.~Verdier,
{\em Construction g\'eom\'etrique de la correspondance de McKay},
Ann. Sci. \'Ecole Norm. Sup. {\bf 16} (1983), 409--449.

\bibitem[Gro]{Gro} I.~Grojnowski,
{\em Instantons and affine algebras I: the Hilbert scheme and
vertex operators}, Math. Res. Lett. {\bf 3} (1996), 275--291.

\bibitem[IN]{IN} Y. Ito and I. Nakamura,
{\em McKay correspondence and Hilbert schemes}, Proc. Japan Acad.
Ser. {\bf  A 72} (1996), 135--138.

\bibitem[LT]{LT} A. Lascoux and J.-Y. Thibon,
{\em Vertex operators and the class algebras of symmetric groups},
Zapiski Nauchnyh Seminarov POMI {\bf 283} (2001), 156--177.

\bibitem[Le]{Le} M. Lehn,
{\em Chern classes of tautological sheaves on Hilbert schemes of
points on surfaces}, Invent. Math. {\bf 136} (1999), 157--207.

\bibitem[LS]{LS} M. Lehn and C. Sorger,
{\em Symmetric groups and the cup product on the cohomology of
Hilbert schemes}, Duke Math. J. {\bf 110} (2001), 345--357.

\bibitem[LQW1]{LQW1} W.-P. Li, Z. Qin and W. Wang,
{\em Ideals of the cohomology rings of Hilbert schemes and their
applications}, Trans. AMS {\bf 356} (2003), 245--265.

\bibitem[LQW2]{LQW2} ---------,
{\em Hilbert schemes, integrable hierarchies, and Gromov-Witten
theory}, Internat. Math. Res. Notices {\bf 40} (2004), 2085--2104.

\bibitem[Mac]{Mac} I.~G. Macdonald,
{\em Symmetric functions and Hall polynomials}, 2nd ed., Clarendon
Press, Oxford, 1995.

\bibitem[McK]{McK} J.~McKay, {\em Graphs, singularities and finite
groups}, Proc. Sympos. Pure Math. {\bf 37}, Amer. Math. Soc,
Providence, RI (1980), 183--186.

\bibitem[MJD]{MJD} T. Miwa, M. Jimbo and E. Date,
{\em Solitons. Differential equations, symmetries and infinite
dimensional algebras}, (originally published in Japanese 1993),
Cambridge University Press, 2000.

\bibitem[Na1]{Na1} H. Nakajima,
{\em Heisenberg algebra and Hilbert schemes of points on
projective surfaces}, Ann. Math. {\bf 145} (1997), 379--388.

\bibitem[Na2]{Na2} ---------,
{\em Lectures on Hilbert schemes of points on surfaces}, Univ.
Lect. Series {\bf 18}, Amer. Math. Soc., 1999.

\bibitem[Na3]{Na3} ---------,
{\em Jack polynomials and Hilbert schemes of points on surfaces},
Preprint, alg-geom/9610021.

\bibitem[OP]{OP} A. Okounkov and R. Pandharipande,
{\em Gromov-Witten theory, Hurwitz numbers, and completed cycles},
math.AG/0204305.

\bibitem[Ru]{Ru} Y.~Ruan,
{\em Stringy geometry and topology of orbifolds}, In: Symposium in
honor of C.~H.~Clemens (Salt Lake City, 2000), Contemp. Math. {\bf
312} (2002), 187--233.

\bibitem[UT]{UT} K. Ueno and K. Takasaki,
{\em Toda lattice hierarchy}, Adv. Studies in Pure Math. {\bf 4}
(1984), 1--95.

\bibitem[Vas]{Vas} E. Vasserot,
{\em Sur l'anneau de cohomologie du sch\'ema de Hilbert de
$\C^2$}, C.~R. Acad. Sci. Paris, S\'er. I Math. {\bf 332} (2001),
7--12.

\bibitem[Wa1]{Wa1} W. Wang,
{\em Algebraic structures behind Hilbert schemes and wreath
products}, In: S.~Berman {\em et al} (eds.), ``Recent developments
in infinite-dimensional Lie algebras and conformal field theory",
Charlottesville, Virginia, May 2000, Contemp. Math. {\bf 297}
(2002), 271--295.

\bibitem[Wa2]{Wa2} ---------,
{\em Vertex algebras and the class algebras of wreath products},
Proc. London Math. Soc. {\bf 88} (2004), 381--404.

\bibitem[Wa3]{Wa3} ---------,
{\em The Farahat-Higman ring of wreath products and Hilbert
schemes}, Adv. in Math. {\bf 187} (2004), 417--446.


\end{thebibliography}
\end{document}